\documentclass[notitlepage, 12pt]{amsart}
\usepackage{amssymb}
\usepackage{amsmath}
\usepackage{amscd}
\usepackage{xspace}
\usepackage[mathscr]{eucal}

\DeclareMathOperator{\sn}{\ensuremath{sn}\xspace}
\DeclareMathOperator{\hhh}{\ensuremath{Hom}\xspace}
\DeclareMathOperator{\eee}{\ensuremath{Ext}\xspace}
\DeclareMathOperator{\enn}{\ensuremath{End}\xspace}
\DeclareMathOperator{\sym}{\ensuremath{Sym}\xspace}
\DeclareMathOperator{\ch}{\ensuremath{ch}\xspace}
\DeclareMathOperator{\CH}{\ensuremath{CH}\xspace}
\DeclareMathOperator{\id}{\ensuremath{id}\xspace}
\DeclareMathOperator{\Hm}{\ensuremath{H}\xspace}
\DeclareMathOperator{\ttt}{\ensuremath{t}\xspace}
\DeclareMathOperator{\D}{\ensuremath{D}\xspace}
\DeclareMathOperator{\hm}{\ensuremath{h}\xspace}

\newcommand{\cal}{\mathcal}
\newcommand{\mapdef}[1]{\ensuremath{\overset{#1}{\longrightarrow}}\xspace}

\newcommand{\bbf}{\mathbb}
\newcommand{\bb}[1]{\ensuremath{\mathbb{#1}}\xspace}

\newcommand{\aaa}{\ensuremath{\alpha}\xspace}

\newcommand{\calc}{\ensuremath{\mathcal{C}}\xspace}

\newcommand{\calo}{\ensuremath{\mathcal{O}}\xspace}

\newcommand{\calp}{\ensuremath{\mathcal{P}}\xspace}

\newcommand{\rar}{\ensuremath{\rightarrow}\xspace}

\newtheorem{theorem}{Theorem}
\newtheorem{prop}{Proposition}
\newtheorem{lemma}{Lemma}
\newtheorem{cor}{Corollary}
\newtheorem{claim}{Claim}

\begin{document}
\title[On the nonexistence of certain morphisms]{On the nonexistence of certain morphisms from Grassmannian to Grassmannian in characteristic $0$ }
\date{\today}
\author{Ajay C. Ramadoss}
\address{Dept. of Mathematics (The University of Chicago)}
\email{ajay@math.uchicago.edu}

\maketitle

\tableofcontents

\section{Introduction}

Problems regarding the constraints that morphisms between homogeneous spaces must satisfy have been studied by Kapil Paranjape and V. Srinivas [7],[8]. In [7], they characterize self maps of finite degree between homogeneous spaces and prove , in addition that finite surjective morphisms from Grassmannian to Grassmannian are actually isomorphisms. In [8] they prove that if $S$ is a smooth quadric hypersurface in ${\mathbb P}^{n+1}$, where $n= 2k+1$, and if $2^k |d$, then there exist continuous maps $f: {\mathbb P}^n \rar S$ so that $f^*(\calo_S(1)) = \calo_{{\mathbb P}^n}(d)$. In the same spirit, one can ask questions like whether there exists a map from a Grassmannian $G(r,n)$ to another Grassmannian $G(r,M)$ so that if $Q$ denotes its universal quotient bundle, $ f^*[Q] = {\psi}^p[Q]$ where $[V]$ denotes the class of a vector bundle $V$ in $K$-theory. Another question in the same spirit would be whether there exist morphisms $f:G(r,n) \rar G(r-1,M)$ so that $f^*(\ch_l(Q)) = ch_l(Q)$. The answers to the first question is in the negative for all $r \geq 2, n \geq 2r+1$ and the answer to the second question is in the negative for infinitely many $r$, with $n$ assumed to be large enough. It may be noted that in these questions, our attention is not restricted solely to dominant/finite morphisms unlike in the results in [7] and [8]. Indeed, the results proven here are not obtainable by the ,methods of [7] and [8] as far I can see. The following theorems proven contain the answers obtained for the above questions. \\

\begin{theorem}
Given any natural number $l \geq 2$,$\exists$ infinitely many natural numbers
$r > 0$ , and a constant $C$ so that if $n > Cr^2+r$, and
$f:G(r,n) \rightarrow G(r-1,\infty)$ is any morphism of schemes, then
$f^*\alpha_l \neq \alpha_l$. Here $\ch([Q]) = \alpha_1 +....+ \alpha_{rank(Q)}+....$.
$Q$ denotes the universal quotient bundle of the grassmannian in consideration
and $\ch$ denotes the chern character. $[Q]$ denotes the class of $Q$ in
$K(G)\otimes \mathbb Q$.
\end{theorem}

\begin{theorem}
$\exists$ infinitely many $r$ so that if $f:G(r,n) \rightarrow G(r-1,\infty)$
is any morphism of schemes with $n > 7r^2 +r+2$, then $f^*\alpha_2 = C{\alpha_1}^2$, $C$ some constant.
\end{theorem}

\begin{theorem}
If $f:G(r,n) \rightarrow G(r,\infty)$ is a morphism of schemes with $r \geq 2$ and
$n \geq 2r+1$, then $f^*[Q] \neq {\psi}^p[Q]$ for any $p \geq 2$
\end{theorem}

In fact, as we prove theorems 1 and 2, we get a little bit more information, which we can state as follows : \\

\begin{theorem}
 Let $F_r\CH^l(G)$ denote the subspace of $\CH^l(G)$ spanned by $\{\ch_l(V)|V a vector bundle of rank \leq r \}$. Then, this filtration is nontrivial as a theory.
In particular, given a fixed $l >0$, there exist infinitely many $r$ such that
$F_{r-1}\CH^l(G(r,n)) \subsetneq F_r\CH^l(G(r,n))$.
\end{theorem}

\begin{theorem}
 If $f:G(3,6) \rar G(2,M) $ is a morphism, then $f^*(\alpha_2(Q')) = C {\alpha_1(Q)}^2$.
\end{theorem}

All these results are proven using certain facts about certain characteristic classes . These characteristic classes were discovered by M. Kapranov[6] (and independently by M.V. Nori [1]) as far as I know. While the proof of these objects being characteristic classes follows [1], the proof of their commuting with Adams operations is provided by me in this paper. These characteristic classes are defined as follows.\\

Let $X$ be a projective variety and let $V$ be a vector bundle on $X$. Consider the Atiyah class $\theta_V \in \Hm^1(X, \enn(V)\otimes \Omega)$. Take the $k$ -fold cup product of $\theta_V$ with itself. Applying the composition map $\enn (V)^{\otimes k} \rightarrow \enn (V)$, followed by the trace map $tr: \enn (V) \rightarrow \calo_X$ to $\theta_V^k$, we obtain an the characteristic class $\ttt_k(V) \in \Hm^k(X, \Omega^{\otimes k})$. Note that the projection $\Omega^{\otimes k} \rightarrow {\wedge}^k \Omega $ when applied to $\ttt_k(V)$ gives us $\ch_k(V)$ which is the degree $k$ part of the Chern character. The classes $\ttt_k$ are referred to in the paper by Kapranov [1] as the {\it big Chern classes}.\\

One may ask about the relation between these classes and the Chern character. The answer to this question lies in the following result presented in this paper.\\

Let $\psi_{X/S}$ is the element of $\psi_{X/S} \in \eee^1(\Omega, \sym^2{\Omega})$ given by the exact sequence  $p_{1*}(0 \rightarrow {\mathcal I}^2/{\mathcal I}^3 \rightarrow {\mathcal I}/{\mathcal I}^3 \rightarrow  {\mathcal I}/ {\mathcal I}^2 \rightarrow 0) =: 0 \rightarrow \sym^2{\Omega} \rightarrow {E} \rightarrow {\Omega} \rightarrow 0$ where ${\mathcal I}$ is the sheaf of ideals defining the diagonal in $X \times_S X$ and $p_1: X \times_S X \rightarrow X$ if the first projection. We also identify $\Hm^k(X, \Omega^{\otimes k})$ with $\eee^k(\calo_X, \Omega^{\otimes k})$. \\

\begin{theorem}
If $X$ is a projective variety and $V$ is a vector bundle on $X$, then $\ttt_k(V) = \ch_k(V) + \sum_{l \leq k} \D_{kl} \circ \ch_l(V)$ where $\D_{kl}$ are elements of $\eee^{k-l}(\Omega^{\otimes l} ,\Omega^{\otimes k})$ which are functorial under pullbacks. Moreover $\D_{kl}$ is obtained from $\psi_{X/S}$ by applying finitely many of the following operations:\\
1: Tensoring with $\id_{\Omega}$\\
2: Yoneda multiplication\\
3: Action of permutation group elements on $\eee^j(-, \Omega^{\otimes j})$.
\end{theorem}

This leads one to check whether the kernels of operations like Yoneda multiplication by $\eee$ elements like $\D_{kl}$ give nontrivial proper subspaces of $\Hm^p(X, \Omega^q)$ in some cases. If that is indeed the case, we get new proper subspaces of the $\Hm^{q,p}$ that are functorial under pullbacks. The following theorem tells us that we indeed have what we are hoping for.   \\

\begin{theorem} For $p,q \geq 2$, let $\Hm^{p,q}_k$ be given by $\ker (\D_{kq} \otimes \id_{\Omega}^{\otimes p-q}) \circ : \Hm^{p,q} \rar  \Hm^k(X, \Omega^{\otimes k+p-q})$ if $p > q$ and  $\ker \D_{kp} \circ : \Hm^{p,q} \rar  \Hm^{k+q-p}(X, \Omega^{\otimes k})$ otherwise. Then, $\Hm^{p,q}_k$ , $k \geq min\{p,q\}$ form an increasing chain of proper subfunctors of $\Hm^{p,q}$
\end{theorem}

At first sight, it may look that theorems 1 and 2 need to be strengthened. Indeed, the on going through the proof, one feels strongly that the filtration, $F_r$ of $CH^l$,which theorem 4 says is nontrivial as a theory, is in fact, strictly increasing as a theory. I feel that   given any $l \geq 2$ fixed, and $r \geq 2$ , there exists some Grassmannian $G=G(r,n)$ so that $\alpha_l(Q) \in F_r \CH^l(G) \setminus F_{r-1} \CH^l(G)$. One approach to this question is entirely combinatorial (along the lines of the proof to theorems 1,2 and 4) and boils down to showing that for some $k$ and a particular $\beta \in KS_k$ depending on $l$ and $k$ only, the subspace spanned by the conjugates of $\beta_{r-1}$ is of strictly smaller dimension than that spanned by conjugates of $\beta_r$. Here, $\beta_i$ is the image of $\beta$ under the projection $KS_k \rar \oplus_{|\lambda| \leq i} End(V_{\lambda})$. Approaching this question along these lines would indeed involve algebraic combinatorics extensively.\\

In a sequel to this paper, I intend to show explicitly how the big Chern classes arise out of the Chern character map to Hochschild homology. This will give us a conceptual way to think about theorem . In addition this will explain the commutativity of these classes with Adams operations in a more conceptual framework. In addition, I intend to use a theorem of Markarian [ ] to show why the $\eee$ elements intrinsic to the variety that arise in the formula for $\ttt_k$ in terms of $\ch_k$ involve $\psi_{X/S}$ in more conceptual terms.\\

 This work would not have been possible without the many useful discussions I had with Prof. M. V. Nori. I am also grateful to Prof. Shrawan Kumar for pointing out a theorem of Bott[4] used in this work and to Prof. Victor Ginzburg for making me aware of the paper by M. Kapranov[6] where the characteristic classes used are introduced. I thank my friend and colleague Apoorva Khare for helping me LaTex this work.   \\

\section{About the Characteristic Classes $\ttt_k$}

Let $V$ be a locally free coherent sheaf on a scheme $X/S$ , with $X$ smooth over $S$. An algebraic connection on $V$ is defined as an $\calo_S$ linear sheaf homomorphism $D:V \rightarrow \Omega_{X/S}\bigotimes_{\calo_X}V$ satisfying the Liebnitz rule, i.e, $D(fv)= df \otimes v +fDv$, where $f \in \Gamma(U,\calo_X)$,$v \in \Gamma(U,V)$, $ \forall U$ open in $X$. Note that a connection on $V$ by itself is not $\calo_X$ linear, but if $D_1$ and $D_2$ are two connections on $V|_U$,$U \subseteq X$ open, then $D_1-D_2 \in \Gamma(U,\enn(V) \otimes \Omega_{X/S})$.\\
For each open $U \subseteq X$ , let $C_V(U)$ denote the set of connections on $V|_U$. This gives us a sheaf of sets on $X$ on which $\enn(V) \otimes_{\calo_X} \Omega_{X/S}$ acts simply transitively. Consider a covering of $X$ by open affines $U_i$, and pick an element $D_i \in C_V(U_i) \forall i$. The $D_i$ together give rise to a well defined element $\theta_V \in \Hm^1(X, \enn(V)\otimes \Omega)$. \\

\begin{lemma} $\theta_{V \otimes W}= A_V+B_W$, where $A_V$ and $B_W$ are the elements in $\Hm^1(X,\enn(V)\otimes \enn(W)\otimes \Omega)$ induced from $\theta_V$ and $\theta_W$ respectively by the maps $\enn(V) \rightarrow \enn(V)\otimes \enn(W)$ ($m \leadsto m \otimes \id_{W}$) and $\enn(W) \rightarrow \enn(V)\otimes \enn(W)$,
($m' \leadsto \id_{V} \otimes m'$) respectively.
\end{lemma}

\begin{cor}
 $\theta_{V \otimes V}$ is induced from $\theta_V$ by the map $\enn(V) \rightarrow \enn(V) \otimes \enn(V)$,($m \leadsto m \otimes \id_{V}+ \id_{V} \otimes m$)
\end{cor}

\begin{proof}
Since $V$ and $W$ are locally free, we can cover $X$ by open sets $U_i$ so that $V$ and $W$ are free over $U_i \forall i$. Let $D_i \in C_V(U_i)$, and $E_i \in C_W(U_i)$ , $\forall i$. The desired result follows from the fact that $\id_V \otimes E_i + D_i \otimes \id_W \in C_(V \otimes W)(U_i)$. \\
\end{proof}

\subsection{Construction of the Characteristic Classes $\ttt_k$}

Given any two locally free coherent sheaves $\cal F$ and $\cal G$ on $X$, one has a cup product $ \cup: \Hm^i(X,\cal F) \otimes \Hm^j(X, \cal G) \rightarrow \Hm^{i+j}(X, \cal F \otimes \cal G)$. Hence, we can consider the cup product of $\theta_V$ with itself $k$ times - $\theta_{V}\cup ......\cup \theta_{V}= \theta_{V}^k \in \Hm^k(X, \enn(V)^{\otimes k} \bigotimes \Omega^{\otimes k})$. The compositon map $ \varphi: \enn(V)^{\otimes k} \rightarrow \enn(V)$ induces a map $\varphi_{*}:\Hm^k(X, \enn(V)^{\otimes k} \otimes \Omega^{\otimes k}) \rightarrow \Hm^k(X, \enn(V) \otimes \Omega^{\otimes k})$. Let ${\ttt_k}\tilde(V) := \varphi_{*}\theta_{V}^k$. Again, the trace map $tr: \enn(V) \rightarrow \calo_{X}$ is $\calo_X$ -linear, and induces $tr_*:\Hm^k(X, \enn(V) \otimes \Omega^{\otimes k}) \rightarrow \Hm^k(X, \Omega^{\otimes k})$. By definition, $\ttt_k(V) := tr_*{\ttt_k}\tilde (V)$. The classes $\ttt_k$ are referred to in Kapranov [6] as the {\it big Chern classes} \\

\subsection{Basic idea of the proofs of theorems 1,2,3 and 4}

Firstly, $\ttt_k$ behaves well with respect to exact sequences. In other words, \\
\begin{lemma} If $0 \rightarrow V' \rightarrow V \rightarrow V'' \rightarrow 0$ is an exact sequence of locally free coherent sheaves on $X$, then $\ttt_k(V) = \ttt_k(V')+\ttt_k(V'')$.
\end{lemma}

This fact shall be proven later.We however, observe that it means that $\ttt_k$ gives rise to a map $\ttt_k:K(X) \otimes \bbf Q \rightarrow \Hm^k(X, \Omega^{\otimes k})$.The following questions arise. \\

Question 1: If $V$ is a vector bundle on $X$ and $\alpha_l(V)$ denotes the projection of $[V] \in K(X) \otimes \bbf Q$ to the eigenspace of $2^l$ under the second adams operator $\Psi^2$ (alternately, $\alpha_l(V)= \ch^{-1}(\ch_l(V))$, where $\ch_l(V)$ is the part of $\ch([V])$ in $\CH^l(X) \otimes \bb Q$). Does $\ttt_k(\aaa_l(V)) \in \Hm^k(X,''\sym^l(L(\Omega))'')$ ?. Here, $L(\Omega)$ is the free super Lie Algebra generated by $\Omega$ in $T\Omega$, which is the universal enveloping algebra of $L(\Omega)$.By the PBW theorem, $L(\Omega)= \bigoplus ``\sym^l(L(\Omega))''$, where $''\sym^l(W)'':= \bigoplus_{p+q=l}{\wedge}^pL_+ \otimes \sym^qL_- $, where $L_+$ and $L_-$ are the even and odd degree parts of $L(\Omega)$ respectively.
It turns out that $\ttt_k(\aaa_l(V))$ lies in the kth cohomology of $X$ with coefficients in a suitable dual of $''\sym^l(L(\Omega))''$. \\

Question 2: Nori[1] shows that the classes $\ttt_k$ are indeed characteristic classes. In other words, they are functorial under pullbacks. Also, $\bb C(S_k)$, the group ring of the permutation group on $k$ letters acts on ${\Omega}^{\otimes k}$ by the natural right action $\sigma:v_1\otimes ....\otimes v_k \leadsto v_{\sigma(1)}\otimes....\otimes v_{\sigma(k)}$.Thus one can talk about the right ideal (in $\bb C(S_k)$ ) of elements $\beta$ such that $\beta_*(\ttt_k(V))=0$. Let $Q_r$ denote the universal qoutient bundle of a Grassmannian $G(r,n)$ (where n is large enough). Let ${\aaa_l}^r := \aaa_l(Q_r)$. Let $I(l,r,k)$ denote the right ideal in $\bb C(S_k)$ annihilating $\ttt_k({\aaa_l}^r)$. Question: For a fixed $l$ and $r$, is there a $k$ so that $I(l,r,k) \subsetneq I(l,r-1,k)$ ?. \\
Note that an affirmative answer to this question implies Theorem 1 (As we shall see later ,theorem 2 is a corollary of Theorem 1). For, suppose there was a morphism $f:G(r,n) \rightarrow G(r-1,M)$ satisfying $f^*({\aaa_l}^{r-1}) = {\aaa_l}^r$, then $\ttt_k(f^*{\aaa_l}^{r-1}) = f^*(\ttt_k({\aaa_l}^{r-1}))= \ttt_k({\aaa_l}^r)$. Also, one notes that if $\beta \in \bb C(S_k)$, then $\beta_*\ttt_k(f^*V)= f^* \beta_*\ttt_k(V)$. Thus, in this situation, $I(l,r-1,k) \subseteq I(l,r,k)$, a contradiction. \\
Note that if the answer to question 2 is in the affirmative, the proof of Theorem 1 can be rephrased as follows: Given any $\beta \in I(l,r-1,k) \setminus I(l,r,k)$, $\beta_*\ttt_k(V)$, is functorial with respect to pullbacks, and $\beta_*\ttt_k({\aaa_l}^r) \neq 0$, and $\beta_*\ttt_k({\aaa_l}^{r-1}) = 0$. Thus, $f^*{\aaa_l}^{r-1} \neq {\aaa_l}^r$. \\
This is the line of thought for proving Theorem 3. We prove the existence of some object in $\Hm^k(X, {\Omega}^{\otimes k})$, for some $k$, obtained by $\ttt_m, m \leq k$ , cup products and the action of the $S_k$ on  $\Hm^k(X, {\Omega}^{\otimes k})$ , that is functorial with respect to pullback, $0$ for $Q$, and nonzero for ${\psi}^pQ$. \\

\subsection{Another question regarding the classes $\ttt_k$}

Note that in general, $\Hm^k(X, \hhh_{\calo}(\aaa,\beta)) = \eee^k(\aaa,\beta)$.For, if $0 \rightarrow \beta \rightarrow I_0 \rightarrow I_1 \rightarrow ......$ is an injective resolution of $\beta$, then $0 \rightarrow \hhh(\aaa, \beta) \rightarrow \hhh(\aaa,I_0) \rightarrow \hhh(\aaa,I_1) \rightarrow .....$ is an injective resolution of $\hhh(\aaa, \beta)$, where $\aaa$ and $\beta$ are locally free coherent sheaves. We can therefore, think of $\tilde{\ttt_k}(V)$ as an element of $\eee^k(V, V \otimes {\Omega}^{\otimes k})$, and $\ttt_k(V)$ as an element of $\eee^k({\calo}_X, {\Omega}^{\otimes k})$. On the other hand, we have a canonical element $\psi_{X/S} \in \eee^1(\Omega, \sym^2{\Omega})$ which is got as follows: Let $\mathcal I$ denote the sheaf of ideals defining the diagonal in $X \times_S X$. Consider the exact sequence $0 \rightarrow {\mathcal I}^2/{\mathcal I}^3 \rightarrow {\mathcal I}/{\mathcal I}^3 \rightarrow  {\mathcal I}/ {\mathcal I}^2 \rightarrow 0$. Consider $p_{1*}(0 \rightarrow {\mathcal I}^2/{\mathcal I}^3 \rightarrow {\mathcal I}/{\mathcal I}^3 \rightarrow  {\mathcal I}/ {\mathcal I}^2 \rightarrow 0) =: 0 \rightarrow \sym^2{\Omega} \rightarrow {E} \rightarrow {\Omega} \rightarrow 0$. This gives us a canonical element in $\eee^1(\Omega, \sym^2{\Omega})$ which we denote by $\psi_{X/S}$. Also, the natural projection $p:{\Omega}^{\otimes k} \rightarrow {\wedge}^k{\Omega}$ induces a map $p_*: \eee^k({\calo}_X, {\Omega}^{\otimes k}) \rightarrow \eee^k({\calo}_X,{\wedge}^k{\Omega} )$, with $p_*(\ttt_k(V)) = \ch_k(V)$. The following question now arises: Is there a formula for $\ttt_k(V)$ in terms of $\ch_l(V) \text{ } \forall  l \leq k$, $\psi_{X/S}$, and permutations applied to $ {\Omega}^{\otimes k} $ ?. The answer to this question is yes. \\

\subsection{Proofs of the basic properties of $\ttt_k$}

In this section we prove the fact that the $\ttt_k$ indeed give us characteristic classes. Besides, we also derive a formula for $\ttt_k(V \otimes W)$ . The results of this section are due to Nori[1]. \\

\begin{lemma}  Let $ 0 \rar V' \rar V \rar V'' \rar 0$ be an exact sequence of locally free coherent sheaves on $X$. Then, $\ttt_k(V) = \ttt_k(V')+ \ttt_k(V'')$.
\end{lemma}

\begin{proof}
Let $V'$,$V$,$V''$, be as in the statement of this lemma. We first prove this lemma for the case when $k=1$. Consider a cover of $X$ by sets $U_i$ such that $V$and $V'$ are trivial over the $U_i$. On each $U_i$, choose a connection $D_i$, so that $D_i|V'$ is a connection on $V'|U_i$ i.e, $D_i(V'|U_i) \subset {\Omega}_{X/S} \otimes V'$. On the other hand, one can consider for each  $U \subseteq X$ open ,\\  $C_{V,V'}(U)= \{\text{Connections on } V|U \text{ that give rise to a connection on } V'|U \}$. Note that ${\mathcal P}\otimes \Omega$ acts simply transitively on  $C_{V,V'}(U)$, where $\mathcal P$ is the subsheaf of sections of $\enn(V)$ that preserve $V'$. Thinking of the $D_i$ as elements of $C_V(U_i)$, we see that they give rise to the Atiyah class $\theta_V  \in \Hm^1(X, \enn(V) \otimes \Omega)$. On the otherhand, when they are thought of as elements of $C_{V,V'}(U_i)$, they give rise to an element $\theta_{V,V'} \in \Hm^1(X, {\mathcal P} \otimes \Omega)$.If $i: {\mathcal P} \rar \enn(V)$ is the natural inclusion, then clearly, ${(i \otimes \id)}_* \theta_{V,V'}= \theta_V $. We shall denote ${i \otimes \id}$ by $i$ henceforth. \\
Note that $tr \circ i = tr $. Hence, $tr_* \theta_{V,V'} = tr_* \theta_V = \ttt_1(V)$. On the other hand, restriction to $V'$ gievs us a map $p_1: \mathcal P \rar \enn(V')$. Then $p_{1*}\theta_{V,V'}$ is the cohomology class obtained by looking at $D_i|V'$ as elements of $C_{V'}(U_i)$ which is $\theta_{V'}$. We also have a projection $p_2: \mathcal P \rar \enn(V'')$. Note that since the $D_i$ are connections on $V$ that restrict to connections on $V'$, they induce connections on $V''$ (all restricted to $U_i$) which we will again denote by $D_i$. Note that $p_{2*}\theta_{V,V'}$ is the cohomology class obtained by thinking of $D_i|{V''}$as elements of $C_{V''}(U_i)$, i.e, $\theta_{V''}$. Now, $tr|_{\mathcal P} = tr \circ p_1 +tr \circ p_2$. This proves the lemma for $k=1$. \\
To prove the lemma in general, one makes the following observations: \\

1. Let $ \tilde{\ttt_k}(V,V') = {\varphi}_* {\theta_{V,V'}}^k$. Then, $i_* \tilde{\ttt_k}(V,V') = \tilde{\ttt_k}(V) $. Where $ \varphi: \enn(V)^{\otimes k} \rightarrow \enn(V)$ is the composition map. This follows from the commutativity of the following diagram: \\

\[
  \begin{CD}
     {\mathcal P}^{\otimes k} @> i^{\otimes k}>> \enn(V)^{\otimes k}\\
      @V{\varphi}VV                                  @VV{\varphi}V\\
     {\mathcal  P}  @> i >>                       \enn(V)
   \end{CD}
  \] \\

2. $p_{1*}\tilde{\ttt_k}(V,V') = \tilde{\ttt_k}(V')$ and $p_{2*}\tilde{\ttt_k}(V,V') = \tilde{\ttt_k}(V'')$. This follows from the commutativity of both  diagrams that follows: \\

\[
  \begin{CD}
\enn(V)^{\otimes k} @> {p_1}^{\otimes k} >> \enn(V')^{\otimes k}\\
    @V{\varphi}VV                               @VV{\varphi}V\\
      \enn(V)       @> p_1 >>                \enn(V')
\end{CD}
\]
   \\

\[
  \begin{CD}
\enn(V)^{\otimes k} @> {p_2}^{\otimes k} >> \enn(V'')^{\otimes k}\\
    @V{\varphi}VV                               @VV{\varphi}V\\
      \enn(V)       @> p_2 >>                \enn(V'')
\end{CD}
\] \\

From this and the additivity of trace, we see that $\ttt_k(V) = \ttt_k(V')+ \ttt_k(V'')$. \\
\end{proof}

\begin{cor} $\ttt_k$ can be extended to a ${\mathbb Q}$ -linear map
 $\ttt_k:K(X) \otimes {\mathbb Q} \rar \Hm^k(X, {\Omega_{X/S}}^{\otimes k})$.
\end{cor}

\begin{lemma}
If $V = V' \oplus V''$ as $\calo_{X}$-modules and $p_1$ and $p_2$ are the natural projections $\enn(V) \rar \enn(V')$ and $\enn(V) \rar \enn(V'')$ respectively, then   $p_{1*}\tilde{\ttt_k}(V) =\tilde{\ttt_k}(V')$ and $p_{2*}\tilde{\ttt_k}(V) = \tilde{\ttt_k}(V'')$.
 \end{lemma}

\begin{proof}
 Pick a cover of $X$ by open sets $U_i$, so that $\exists D'_i \in C_{V'}(U_i)$ and $D''_i \in C_{V''}(U_i)$. Then $D_i = D'_i \oplus D''_i \in C_V(U_i)$. The rest of the proof is got by following the proof of Lemma 2. \\
\end{proof}

Another fact that needs to be checked is that $\ttt_k$ is fuctorial under pullback. If $f:X \rar Y$ is a morphism, then we have the pullback $f^*:K(Y) \rar K(X)$. \\
Also, we have a map $ \eta: f^*{\Omega_Y} \rar {\Omega_X}$ and a map $ f^*: \Hm^k(X, {\Omega_Y}^{\otimes k}) \rar \Hm^k(X, {f^*{\Omega_Y}}^{\otimes k})$. Define, by abuse of notation, $f^*: \Hm^k(Y, {\Omega_Y}^{\otimes k}) \rar \Hm^k(X, {\Omega}^{\otimes k})$ to be $\eta \circ f^* $. We now have: \\

\begin{lemma}
$\ttt_k(f^*V) =f^*(\ttt_k(V))$
\end{lemma}

\begin{proof}
We first show this result for $k=1$. Suppose $Y$ is covered by open sets $U_i$ and $D_i \in C_V(U_i)$. Then, on $f^{-1}U_i$, the pull-back connection $f^*D_i$ is the canonical connection that takes a section $f^.v$ to $\eta(f^.(D_i(v)))$. Thus, $\theta_{f^*V} = \eta_*f^*\theta_{V}$. This shows the result for $k=1$. The rest follows from the commutativity of the following diagram:

\[
 \begin{CD}
 \enn(V)^{\otimes k} @> {f^*}^{\otimes k} >> \enn(f^*V)^{\otimes k}\\
 @V{\varphi}VV                                   @VV{\varphi}V\\
  \enn(V)            @> f^* >>               \enn(f^*V)\\
  @VtrVV                                        @VVtrV\\
 \calo_{X}          @> f^* >>               \calo_{Y}
\end{CD}
\]
\\
\end{proof}

Having done this, we still need to show the alternate description of $\tilde{\ttt_k}(V)$ and $\ttt_k(V)$ as elements of $\eee^k(V, V \otimes {\Omega}^{\otimes k})$ and $\eee^k({\calo_{X}},{\Omega}^{\otimes k})$ respectively. This proof we shall relegate to the appendix. Another important property that we prove here is that $ \oplus \ttt_k: K(X)\otimes {\mathbb Q} \rar \oplus \Hm^k(X, {\Omega}^{\otimes k})$ is a ring homomorphism provided that the right hand side is equipped with a suitable multiplicative structure that we shall describe explicitly. \\

\begin{lemma}
If $V$ and $W$ are two locally free coherent sheaves on $X$, then,
$\ttt_k(V \otimes W) = \sum_{l+m=k} \ttt_l(V) \odot \ttt_m(W)$, where $\odot$ is the product $\Hm^l(X, {\Omega}^{\otimes l}) \otimes \Hm^m(X, {\Omega}^{\otimes m}) \rar \Hm^k(X, {\Omega}^{\otimes k})$ induced by the signed shuffle product ${\Omega}^{\otimes l} \otimes  {\Omega}^{\otimes m} \rar  {\Omega}^{\otimes l+m}$.
\end{lemma}

\begin{proof}
 We know that $\theta_{V \otimes W} = \theta_{V} \otimes \id_W + \id_V \otimes \theta_{W}$. Therefore, ${\theta_{V \otimes W}}^k = (A_V + B_W) \cup .... \cup (A_V+B_W) $, where $A_V = \theta_{V} \otimes \id_W$ and $B_W = \id_V \otimes \theta_{W}$.  Thus, ${\theta_{V \otimes W}}^k = {(A_V+B_W)}^k = \sum_{l+m =k} \sum_{\sigma \text{ a (l,m)-shuffle }} \sn (\sigma){{\sigma}^{-1}}_* {A_V}^l \cup {B_W}^m$. \\
This requires some elaboration. Here, an $(l,m)$-shuffle is a permutation of ${\{1,2,3,....,l+m\}}$ such that $ \sigma(1) < ....< \sigma(l)$ and $\sigma(l+1)<....<\sigma(l+m)$. Also, a given permutation $\mu \in S_k$ acts on $\enn(V \otimes W)^{\otimes k} \bigotimes {\Omega}^{\otimes k}$ by $v_1 \otimes ....\otimes v_k \bigotimes w_1 \otimes .....\otimes w_k \leadsto v_{\mu(1)} \otimes..... \otimes  v_{\mu(k)} \bigotimes  w_{\mu(1)} \otimes.....\otimes  w_{\mu(k)}$ and therefore induces a map from $\Hm^k(X,{\enn(V \otimes W)}^{\otimes k} \otimes {\Omega}^{\otimes k})$ to itself. Now, we can see that in $(A_V+B_W)^k$, terms having $l$ $A_V$'s cupped with $m$ $B_W$'s are in one-one correspondence with sequences $b_1 <....<b_m , b_i \in \{1,2,3,...,l+m\} \text{ } \forall i $ (The $b_i$'s being the positions of the $B_W$'s) which are in turn in one-one correspondence with $(l,m)$ shuffles. Now, it is easy to see that if $ \sigma$ is the $(l,m)$-shuffle given by $\sigma(l+i) = b_i , 1 \leq i \leq m$, then, applying $ \sigma_*$ to the term with the $B_W$'s placed in the positions $b_i$ gives us $\sn(\sigma){A_V}^l \cup {B_W}^m$. The lemma is now proven by recognising that $tr_* \circ {\varphi}_* {\sigma}_* ({A_V}^l \cup {B_W}^m) = {\sigma}_* \ttt_l(V) \cup \ttt_m(W)$ if $\sigma$ is a $(l,m)$-shuffle. This is because an $(l,m)$-shuffle does not change the order of composition among the $\enn(V)$-terms and among the $\enn(W)$ terms respectively. \\
\end{proof}

Note that the signed shuffle product on $T{\Omega}$, the tensor co-algebra on ${\Omega}$, makes  $T{\Omega}$ into a commutative Hopf-algebra. Therefore, $\odot$ gives a commutative ring structure on $\oplus \Hm^k(X, {\Omega}^{\otimes k})$, and Lemma 6 tells us that  $ \oplus \ttt_k: K(X)\otimes {\mathbb Q} \rar \oplus \Hm^k(X, {\Omega}^{\otimes k})$ is a ring  homomorphism. \\

\section{$\lambda$-ring structure on $\oplus \Hm^k(X, {\Omega}^{\otimes k})$ and the commuting of $\ttt_k$ with respect to Adams operations}

Here, we show that $\oplus \Hm^k(X, {\Omega}^{\otimes k})$ has a special $\lambda$-ring structure (i.e, has Adams operations) and prove that the classes $\ttt_k$ commute with Adams operations. This will basically answer Question 1 of the previous section in the affirmative . Before proceeding, we need a digression on Hopf-algebras. \\

\subsection{Adams operations on commutative Hopf-algebras}

This material is from Loday[2]. If ${\mathcal H} =({\mathcal H},\mu,\Delta,u,c)$ is a commutative Hopf-algebra over a field $K$, we can define the convultion of two maps $f,g \in \enn_K({\mathcal H})$ by $f*g = \mu \circ (f \otimes g) \circ  \Delta$. The convolution product $*$ is an associative product on $\enn_K({\mathcal H})$. One can therefore consider operations ${\psi}^k := \id*...*\id \in \enn_K({\mathcal H})$, i.e, the identity convolved with itself $k$ times. These are $K$-linear, and satisfy ${\psi}^p \circ {\psi}^q = {\psi}^{pq}$. If they turn out to be ring homomorphisms as well, then they can be thought of as Adams operations. For this, we need to use the fact that ${\mathcal H}$ is a commutative Hopf-algebra.  In this situation, one checks that if $f:{\mathcal H} \rar {\mathcal H}$ is a ring homomorphism, so is $\id*f$. Also observe that the map $uc$ is the identity with respect to $*$. If ${\mathcal H}$ is graded, with ${\mathcal H}_0 = K$, then $K$-linear maps of degree $0$ form a subalgebra of $(\enn_K({\mathcal H}), +,*)$. If $f:{\mathcal H} \rar {\mathcal H}$,$f(1) =0 \implies f^{*k}|_{{\mathcal H}_n} = 0 \text{ } \forall n<k$. Hence the series $e^{(1)}(f) :=  \ln (uc+f) = f-f^2/2 +f^3/3 -......$ is a polynomial when restricted to each graded part of ${\mathcal H}$, and thus makes sense. For the same reason, $e^{(i)}(f):= {e^{(1)}(f)}^i/i!$ makes sense. For each $n$, this determines a $K$-linear endomorphism ${e^{(i)}}_n(f)$ of ${\mathcal H}_n$. Clearly, ${e^{(1)}}_0(f) =0$ and ${e^{(i)}}_n(f)=0 \forall i >n$. Since $(1+x)^k= \exp(k \ln (1+x))$, we get $(uc +f)^k = uc + \sum_{i \geq 1} k^ie^{(i)}(f)$. Putting $f= \id -uc$ and observing that $f(1)=0$, we see that ${\psi}^k = \sum_{i \geq 1}k^i{e^{(i)}}_n$, where ${e^{(i)}}_n = {e^{(i)}}_n(f)$.It is easily seen that the elements ${e^{(i)}}_n$ are mutually orthogonal idempotents adding up to $\id$ in $\enn_K({\mathcal H}_n)$. \\
 The Hopf algebra that is relevant to us is the tensor co-algebra of a vector space $V$. Here, $T^*(V)_n = V^{\otimes n}$, $\Delta(v_1 \otimes .....\otimes v_n) = \sum_{0 \leq i \leq n} v_1 \otimes .....\otimes v_i \bigotimes v_{i+1} \otimes .....\otimes v_n $ (cut coproduct),and $ \mu(v_1 \otimes .....\otimes v_p \bigotimes v_{p+1} \otimes .....\otimes v_{p+q}) = \sum_{\sigma \text{ a (p,q)-shuffle}} v_{{\sigma}^{-1}(1)} \otimes .....\otimes v_{{\sigma}^{-1}(p)} \otimes v_{{\sigma}^{-1}(p+1)} \otimes .....\otimes v_{{\sigma}^{-1}(p+q)}$ We therefore, note that in this case, ${\psi}^2(v_1 \otimes .....\otimes v_n) = \sum_{p+q=n}\sum_{\sigma \text{ a (p,q)-shuffle}} v_{{\sigma}^{-1}(1)} \otimes .....\otimes v_{{\sigma}^{-1}(p)} \otimes v_{{\sigma}^{-1}(p+1)} \otimes .....\otimes v_{{\sigma}^{-1}(n)}$.In this particular case, we also want to find out about the idempotents ${e^{(i)}}_n \in \enn_K(V)^{\otimes n}$. The following result from Loday[2] is exactly what we want. \\

\begin{lemma}
$e^{(i)}_n = \sum_{j=1}^n {a_n}^{i,j}{l_n}^j $ where $ \sum_{i=1}^n {a_n}^{i,j}X^i = \binom{X-j+n}{n}$ and ${l_n}^j = \sum_{\sigma \in S_{n,j}} (\sn \sigma) {{\sigma}_{*}}^{-1}$. Here, $S_{n,j} =\{ \sigma \in S_n | card\{i| \sigma(i) > \sigma(i+1) \} = j-1 \}$, ${{\sigma}_{*}}^{-1}( v_1 \otimes ..... .....\otimes v_n) =  v_{{\sigma}^{-1}(1)} \otimes ..... .....\otimes v_{{\sigma}^{-1}(n)}$ .
\end{lemma}

For example,${e^{(n)}_n} = \sum_{\sigma \in S_n} \sn (\sigma) {{\sigma}_{*}}^{-1}$.\\

\subsection{Description of $\lambda$-ring structure on $\oplus \Hm^k(X, {\Omega}^{\otimes k})$}

Consider the tensor co-algebra $T^*\Omega$. Consider the Adams operations ${\psi}^k$ on $T^*\Omega$ as described in the previous subsection. Note that ${\psi}^k|_{{\Omega}^{\otimes n}}$ induces a map ${\psi}^k_* : \Hm^n(X, {\Omega}^{\otimes n}) \rar  \Hm^n(X, {\Omega}^{\otimes n}) $. Thus the Adams operation ${\psi}^k$ induces a map ${\psi}^k_*: \oplus \Hm^n(X, {\Omega}^{\otimes n}) \rar \oplus \Hm^n(X, {\Omega}^{\otimes n})$ that is $K$-linear. That ${\psi}^p \circ {\psi}^q = {\psi}^{pq}$ implies that $ {\psi}^p_* \circ {\psi}^q_* = {\psi}^{pq}_*$. Define the $k$-th Adams operation on $\oplus \Hm^n(X, {\Omega}^{\otimes n})$ to be ${\psi}^k_*$. That the Adams operations so defined commute with multiplication in this ring follows from Lemma 6 of Section 2, which says that multiplication in this ring is induced by the multiplication in $T^*\Omega$. We have therefore, proven the following Lemma: \\

\begin{lemma}
$\oplus \Hm^n(X, {\Omega}^{\otimes n})$ is a special $\lambda$-ring with Adams operations ${\psi}^p$ given by ${\psi}^p_*$.
\end{lemma}

We now begin proving the fact that the classes $\ttt_k$ commute with Adams operations. By the corollary to Lemma 1, ${\theta}_{V \otimes V}$ is induced from ${\theta}_V$ by the map $\beta:\enn(V) \rar \enn(V)$ given by $m \rar m \otimes \id_V + \id_V \otimes m$ i.e, ${\theta}_{V \otimes V} = {\beta}_*{\theta}_V$. Therefore, ${{\theta}_{V \otimes V}}^k = {\beta}_*{\theta}_V \cup .....\cup {\beta}_*{\theta}_V = ( \beta \otimes ....\otimes \beta)_* {\theta_V}^k$. By abuse of notation, we shall refer to $\beta \otimes ....\otimes \beta$ as $\beta$. Then, ${\theta}_{V \otimes V} = {\beta}_* {\theta}_V$, where $\beta: \enn(V)^{\otimes k}  \rar \enn(V)^{\otimes k}$ is given by $m_1 \otimes m_k \rar {\bigotimes_{i=1}}^k (m_i \otimes \id_V +\id_V \otimes m_i)$. Further, a direct computation shows that if $W$ is a vector space over a field $F$, with $char F \neq 2$, $W \otimes W = \sym^2 W \oplus {\wedge}^2 W$ and if $p_1$ and $p_2$ are the resulting projections from $\enn(W) \otimes \enn(W) = \enn(W \otimes W)$ onto $\enn(\sym^2W)$ and $\enn({\wedge}^2W)$ respectively, if $M , N \in \enn(W)$, then $tr(p_1( M\otimes N)) - tr(p_2( M \otimes N)) = tr(M \circ N)$. By this fact, and Lemma 4, we see that $\ttt_k( {\psi}^2V) = \ttt_k(\sym^2 V) - \ttt_k({\wedge}^2 V) = tr_* p_{1*} \tilde{\ttt_k}(V \otimes V) -tr_*p_{2*} \tilde{\ttt_k}(V \otimes V) = tr_* {\alpha}_* \tilde{\ttt_k}(V \otimes V)$ , where $\alpha : \enn(V) \otimes \enn(V) \rar \enn(V)$ is the composition map. \\
Let $ {\varphi}: \enn(V \otimes V)^{\otimes k} \rar \enn(V \otimes V)$ be the composition map. Observe that $ \alpha \circ \varphi \circ \beta : \enn(V)^{\otimes k} \rar \enn(V)$ is the map given by $m_1 \otimes ... \otimes m_k \leadsto  \sum_{p+q =k} \sum_{ \sigma \text{ a (p,q) -shuffle} } m_{\sigma(1)} \circ ... \circ m_{\sigma(k)}$ ($\circ$ denoting the usual matrix multiplication on the right hand side of the last equation). Consider the map $ \gamma:\enn(V)^{\otimes k} \rar \enn(V)^{\otimes k}$ given by $ m_1 \otimes....\otimes m_k \leadsto  \sum_{p+q =k} \sum_{ \sigma \text{ a (p,q) -shuffle} } m_{\sigma(1)} \otimes ... \otimes m_{\sigma(k)}$. Then, we see that $tr_* \circ \varphi_* \circ \gamma_* {\theta_V}^k = tr_* \circ \alpha_* \tilde{\ttt_k}(V \otimes V) = \ttt_k({\psi}^2V)$. Also observe that ${\psi}^2\ttt_k(V) = tr_* \varphi_* {\psi}^2_* {\theta_V}^k$ since the following diagram commutes: \\

\[
 \begin{CD}
   \enn(V)^{\otimes k} \otimes {\Omega}^{\otimes k}  @> \id \otimes {\psi}^2 >>    \enn(V)^{\otimes k} \otimes {\Omega}^{\otimes k} \\
    @V{tr \circ(\varphi \otimes \id)}VV                                            @VV{tr \circ (\varphi \otimes \id)}V \\
      {\Omega}^{\otimes k}                 @>{\psi}^2>>                            {\Omega}^{\otimes k}
    \end{CD}
\] \\

Here ${\psi}^2_*$ on $\Hm^k(X, \enn(V)^{\otimes k} \otimes {\Omega}^{\otimes k})$ is by definition induced by the ${\psi}^2$ on $ {\Omega}^{\otimes k}$. Thus, the following lemma remains to be proven: \\

\begin{lemma}
   $\gamma_* {\theta_V}^k = {\psi}^2_* {\theta_V}^k$
\end{lemma}

\begin{proof}

Note that the cup-product is anti-commutative. Therefore, if $\sigma \in S_k$ , then the map given by $\sigma: m_1 \otimes....\otimes m_k \bigotimes v_1 \otimes......\otimes v_k \leadsto  \sn(\sigma) m_{\sigma(1)} \otimes ... \otimes m_{\sigma(k)} \bigotimes  v_{\sigma(1)} \otimes ... \otimes v_{\sigma(k)}$ preserves ${\theta_V}^k$. It now suffices to note that $ \gamma = {\psi}^2 \circ \sigma$. \\
\end{proof}

With this, we have proven the following lemma: \\

\begin{lemma}
$\ttt_k({\psi}^2V)= {\psi}^2\ttt_k(V)$.
\end{lemma}

Recalling that $\alpha_l(V) = \ch^{-1}(\ch_l(V))$, where $\ch$ is the Chern character map, we now have the corollary below:\\
\begin{cor}
$\ttt_k(\alpha_l(V)) ={e_k}^{(l)}_*\ttt_k(V)$ where ${e_k}^{(l)}$ is the idempotent described in Lemma 6.
\end{cor}

\begin{proof}
Note that ${\psi}^2 = \sum e^{(l)}2^l$. The fact that the ${e_k}^{(l)}$ are mutually orthogonal idempotents adding upto $\id$ tells us that ${\psi}^2 \circ {e_k}^{(l)} = 2^l {e_k}^{(l)}$. Therefore, ${\psi}^2 \ttt_k(V) = \sum  2^l {e_k}^{(l)}_* \ttt_k(V) = \ttt_k({\psi}^2 V) = \ttt_k( \sum 2^l \alpha_l(V)) = \sum 2^l \ttt_k(\alpha_l(V))$. Since eigenvectors corresponding to different eigenvalues of a linear operator on a finite dimensional vector space over a field of characteristic $0$ are linearly independent, the desired result follows.\\
\end{proof}

Remark: The lemma stating that the classes $\ttt_k$ commute with Adam's operations answers Question 1. of Section 2. To be more precise, if $TV$ is the graded tensor algebra over a vector space $V$, (with usual tensor product giving the multiplication, and coproduct dictated by the fact that $V \subset TV$ are primitive elements), then $T^*V$ is the graded Hopf algebra dual to $TV$. The map ${\psi}^2 = \mu \circ \Delta: T^*V \rar T^*V$ has as its dual the map $ \mu \circ \Delta:TV \rar TV $. The $2^l$-eigenspace of this map is seen to be $''\sym^l(L(V))''$. Thus, the $2^l$-eigenspace of ${\psi}^2:T^*V \rar T^*V$ is dual to the space $''\sym^l(L(V))''$. Thus, $\ttt_k(\alpha_L(V))$ lands in $k$-cohomology with coefficients in a space dual to $''\sym^l(L(\Omega))''$. Moreover, the last corollary explicitly describes the projector that gives $\ttt_k(\alpha_l(V))$ from $\ttt_k(V)$ as the action on $\ttt_k(V)$ of a certain idempotent in ${\mathbb C}(S_k)$.\\

\section{Calculating $\ttt_k(Q)$, $Q$ the universal quotient bundle of a Grassmannian $G(r,n)$}

 In this section, we compute $\ttt_k(Q)$ and $\tilde{\ttt_k}(Q)$ as elements in $\eee^k(K, {\Omega}^{\otimes k})$ and  $\eee^k(Q, Q \otimes {\Omega}^{\otimes k})$ respectively. Here $\eee$ means $\eee$ in the category of $P$-representations, where $G(r,n) = Gl_K(n)/P$, $P$ a parabolic subgroup. This is good enough because all the cohomology classes that we are computing are elements in $\eee^k(-,-)$  (where $\eee$ is in the category of locally free coherent sheaves on $Gl_K(n)/P$) given by exact sequences of $G$ -equivariant vector bundles, where $G=Gl_K(n)$. It follows that the classes $\ttt_k(Q)$ are elements of ${\eee^k(\calo_{G/P},{\Omega}^{\otimes k})}^G$ .It is a theorem of Bott [4] that ${\eee^k(\calo_{G/P},\bar{\calo_G \otimes V})}^G$ is isomorphic to $\eee^k_{P \text{ representations}}(K,V)$ where $K$ is the base field. Here, $V$ is any $P-$ representation, and, for any sheaf $\mathcal F$ with descent data on $G$, $\bar{\mathcal F}$ is the sheaf given by the presheaf $\bar{\mathcal F}(U) = \{ s \in \Gamma(f^{-1}(U),\mathcal F) | I(p_1,p_2) p_1^*s = p_2^*s\}$ where $p_1$ and $p_2$ are the natural projections $G \times_{G/P} G$ and $I(p_1,p_2): p_1^* \mathcal F \rar p_2^* \mathcal F$ is the isomorphism given in the descent data. To do this, we need to use alternate definitions on $\ttt_k(V)$ and $\tilde{\ttt_k}(V)$ as elements in $\eee^k(\calo_X, {\Omega}^{\otimes k})$ and  $\eee^k(V, V \otimes {\Omega}^{\otimes k})$  respectively. \\

\subsection{Alternate construction for $\tilde{\ttt_k}(V)$ and $\ttt_k(V)$}

Let $V$ be a locally free coherent sheaf on a (separated) scheme $X/S$. Then, we have a diagonal embedding $\Delta: X \rar Y:= X \times_S X$. Let $\mathcal I$ be the sheaf of ideals defining this diagonal. Let $p_1$ and $p_2$ denote the two canonical projections from $Y$ to $X$. We have an exact sequence $0 \rightarrow {\mathcal I}/{\mathcal I}^2 \rightarrow {\calo_Y}/{\mathcal I}^2 \rightarrow  {\calo_Y}/ {\mathcal I} \rightarrow 0$. Consider the exact sequence $0 \rightarrow {p_2}^*V \otimes {\mathcal I}/{\mathcal I}^2 \rightarrow {p_2}^*V \otimes {\calo_Y}/{\mathcal I}^2 \rightarrow  {p_2}^*V \otimes {\calo_Y}/ {\mathcal I} \rightarrow 0$ ($\otimes$ denoting tensoring over $\calo_Y$) of coherent $\calo_Y$-modules. Since $ {p_2}^*V \otimes {\mathcal I}/{\mathcal I}^2 $ is supported on the diagonal only, and $p_1$ induces an isomorphism on the diagonal, $R^ip_{1*}( {p_2}^*V \otimes {\mathcal I}/{\mathcal I}^2 ) =0 \text{ } \forall i >0$. Therefore, applying $p_{1*}$, we get an exact sequence  $0 \rightarrow p_{1*}({p_2}^*V \otimes {\mathcal I}/{\mathcal I}^2) \rightarrow p_{1*}({p_2}^*V \otimes {\calo_Y}/{\mathcal I}^2) \rightarrow  p_{1*}({p_2}^*V \otimes {\calo_Y}/ {\mathcal I}) \rightarrow 0 $ of coherent $\calo_X$ -modules. This is an exact sequence $0 \rar V \otimes \Omega \rar J_1(V) \rar V \rar 0$. This exact sequence gives rise to an element in $\eee^1(V, V \otimes \Omega) = \Hm^1(X, \enn(V) \otimes \Omega)$. To see that this element is indeed ${\theta_V}$, note that splittings of the exact sequence  $0 \rar V \otimes \Omega \rar J_1(V) \rar V \rar 0$ on an open subset $U \subset X$ correspond to connections on $U$. Therefore, if $U_i$ form an open cover of $X$, with $D_i$ a connection on $V|_{U_i}$, each $D_i$ corresponds to a splitting of $0 \rar V \otimes \Omega \rar J_1(V) \rar V \rar 0$ restricted to $U_i$. Call this splitting map $f_i: V \rar J_1(V)$. Clearly, on $U_i \cap U_j$,  $f_i-f_j = D_i-D_j : V \rar V \otimes \Omega$. We now only need to use the fact that if $0 \rar \mathcal A \rar \mathcal B \rar \mathcal C \rar0$ is an exact sequence of locally free coherent sheaves on $X$, and $f_i: \mathcal C \rar \mathcal B$ is a splitting of this exact sequence restricted to $U_i$ for each $i$, then the element of $\Hm^1(X,\hhh(\mathcal C,\mathcal A))$ given by the Cech cocycle $(f_i-f_j)$ is indeed the element of $\eee^1(\mathcal C, \mathcal A)$ given by $0 \rar \mathcal A \rar \mathcal B \rar \mathcal C \rar0$. We also observe that if $\alpha \in \Hm^i(X, \mathcal F) = \eee^i(\calo_X , \mathcal F)$ is given by an exact sequence $0 \rar \mathcal F \rar Y_1 \rar ...\rar Y_i \rar \calo_X \rar 0$ and if $\beta \in \Hm^j(X, \mathcal G) = \eee^j(\calo_X , \mathcal G)$ is given by an exact sequence $0 \rar \mathcal G \rar Z_1 \rar ...\rar Z_j \rar \calo_X \rar 0$, then the product $*$ such that $\alpha * \beta$ is the element in $\Hm^{i+j}(X, \mathcal F \otimes \mathcal G) = \eee^{i+j}(\calo_X, \mathcal F \otimes \mathcal G)$ defined by the exact sequence which is the tensor product of the exact sequences representing $\alpha$ and $\beta$ respectively, has the linearity and anticommutativity properties required of the cup product. Since all the cohomology classes we are dealing with are of this type, we can define the cup product to be the product $*$.With this definition of the cup product, it will follow that $\tilde{\ttt_k}(V) \in \eee^k(V, V \otimes {\Omega}^{\otimes k})$ is given by $({\theta_V} \otimes {\id_{\Omega}}^{k-1}) \circ .... \circ {\theta_V}$ where $\circ$ denotes the Yoneda product and ${\theta_V}$ is treated as an element in $\eee^1(V, V \otimes \Omega)$. \\

\subsection{Computation of $\tilde{\ttt_1}(Q)$}

Next, we need to show that in order to calculate $\ttt_k(Q)$ and $\tilde{\ttt_k}(Q)$ as elements in $\eee^k(K, {\Omega}^{\otimes k})$ and  $\eee^k(Q, Q \otimes {\Omega}^{\otimes k})$ respectively (in the category of $P$-representations) it is enough to perform our calculation in the category of $N$-representations, where $N$ is the unipotent subgroup associated with $P$. This is because we have a first quadrant spectral sequence (Lyndon-Hochschild-Serre spectral sequence) \\
$$ {E_2}^{pq} =\Hm^p(P/N ;\Hm^q(N ;A)) \implies \Hm^{p+q}(P;A)$$
for any $P$-representation $A$. But $P/N = Gl(Q) \times Gl(S)$ and all finite-dimensional $Gl(Q) \times Gl(S)$- representations are semisimple. Therefore, for any $P/N$-representation $B$, $\hhh(B,-)$ is exact. Thus, the Lyndon-Hochschild-Serre spectral sequence collapses, leaving only the bottom row. This implies that ${({\eee^*}_{N - representations})}^{P/N} = {\eee^*}_{P-representations}$. This tells us that it suffices to perform our calculations in the category of $N$-representations.  We know that $Q$,$S$ and hence $\Omega$ are all $P$ -representations. We must now see that $\theta_Q$ can indeed be represented in $\eee^1(Q, Q \otimes \Omega)$ by an exact sequence of $P$-representations with $P$-module homomorphisms. However, we note that in the category of locally free coherent sheaves on $G=G(r,n)$, the folowing diagram commutes. \\

  \[
   \begin{CD}
       0 @> >> S @ > \varphi >>    V   @ > \epsilon >> Q @ > >> 0 \\
      @.         @V{\Delta}VV       @VVV                 @VV{\id}V \\
       0 @> >> Q \otimes Q^* \otimes S @> >> J_1(Q) @ > \gamma >> Q @> >> 0
   \end{CD}
\] \\

Here, $ \Delta^* := ev \otimes \id: Q^* \otimes Q \otimes S^* \rar S^*$ , $ev: Q^* \otimes Q \rar K$ being the natural eveluation map. The top row of this diagram is the exact sequence giving $\theta_V$. By the universal property of push-forwards, we see that the following diagram commutes : ($F$ denotes the pushforward $V \amalg_S Q^* \otimes Q \otimes S$).\\

\[
 \begin{CD}
   0 @> >> S @ > \varphi >>  V @> \epsilon >> Q @ > >> 0\\
   @.   @V{\Delta}VV    @VVV             @VV{\id}V\\
  0 @> >> Q^* \otimes Q \otimes S @> >> F @ > >> Q @> >> 0\\
   @.    @VV{\id}V                       @VVV       @V{\id}VV\\
    0 @> >> Q^* \otimes Q \otimes S @> >> J_1(Q) @ > >> Q @> >> 0
\end{CD}
\] \\

Therefore, $\theta_Q$ can be represented by the second row of the above diagram in $\eee^1(Q, Q \otimes \Omega)$. Observe, however, that every arrow in this exact sequence is a $P$-module homomorphism (of course, $Q^* \otimes Q \otimes S$, $V$ and therefore $F$ are all $P$-modules). Thus $\theta_Q$ can be represented by an exact sequence in the category of $P$-representations. It follows that for all $k \geq 1$ , $\tilde{\ttt_k}(Q)$ and $\ttt_k(Q)$ can be represented by  exact sequences in the category of $P$-representations. We next note that $N$ is a Lie group, and the category of $N$-representations is equivalent to the category of $\eta$-representations, where $\eta$ is the Lie-algebra associated to $N$. Also, the category of $\eta$-representations is equivalent to the category of $U(\eta)$-representations, where $U(\eta)$ is the universal enveloping algebra of $\eta$. Since $\eta$ is abelian, $U(\eta) = \sym^*\eta$. In what follows, we shall work in the category of $\sym^*\eta$-modules. Note that as a vector space, $\eta = \Omega =Q^* \otimes S$. Since $N$ acts trivially on $Q$, the ideal $ \Omega \oplus \sym^2\Omega \oplus ....$ annihilates $Q$ thought of as a $\sym^*\Omega$-module. Therefore, a projective resolution of $Q$ can be obtained by taking the Koszul complex $ ....\rar Q \otimes {\wedge}^k\Omega \otimes \sym^*\Omega \rar Q \otimes {\wedge}^{k-1}\Omega \otimes \sym^*\Omega \rar .... \rar Q \otimes \sym^* \Omega \rar Q \rar 0$. It follows that if $W$ is any other $\sym^*\Omega$-module, then $\eee^k(Q,W)$ is just the $k$-th cohomology of the complex $0 \rar \hhh( Q \otimes \sym^* \Omega ,W) \rar ... \rar... \hhh(Q \otimes {\wedge}^K\Omega \otimes \sym^* \Omega,W) \rar...$. If $W$ is also a trivial $\sym^*\Omega$-module, then we see that $\hhh(Q \otimes {\wedge}^K\Omega \otimes \sym^* \Omega,W) = \hhh_K(Q \otimes {\wedge}^k\Omega,W)$ and the Koszul differential in the previous complex is $0$. In this case, $\eee^k(Q,W) = \hhh_K(Q \otimes {\wedge}^k\Omega,W)$. Therefore, to find $\theta_Q$, we need to find arrows $\alpha$ and $\beta$ so that all squares in the following diagram commute: \\

\[
 \begin{CD}
    0 @> >> Q^* \otimes Q \otimes S @> >> F @ > >> Q @> >> 0\\
     @.    @A{\Delta}AA      @AAA             @AA{\id}A     \\
    0 @> >> S @ > \varphi >>  V @> \epsilon >> Q @ > >> 0\\
    @.        @A{\alpha}AA     @A{\beta}AA      @AA{\id}A     \\
    ...  @> >> Q \otimes \Omega \otimes \sym^*\Omega @> >>Q \otimes \sym^*\Omega @> >> Q @> >> 0
\end{CD}
\] \\

Observe that $\Omega = \hhh_K(Q,S) \subseteq \enn(V)$ (Here, we have chosen a $K$-vector space splitting $0 \rar S  \rar V \leftrightarrows Q \rar 0$. Choosing such a splitting describes $\Omega$ as the subspace of elements in $\enn(V)$ consisting of matrices whose ``upper right block'' is the only nonzero block. Note that the product of two such matrices is $0$. Thus, any element of $\sym^*\Omega$ can be thought of as an elemnent of $\hhh(Q,V) \subset \enn(V)$. In this scheme of things, we choose $\beta$ to be the natural eveluation map, and $\alpha$ the restriction of $\beta$ to $Q \otimes \Omega \otimes \sym^*\Omega$. Note that $\beta$ and $\alpha$ are $\sym^*\Omega$-module homomorphisms by construction. Note that $\alpha: Q \otimes \Omega \otimes \sym^*\Omega$ is the $\sym^* \Omega$-module homomorphism induced by $\tilde{\alpha}:= ev \in \hhh_K(Q \otimes\Omega, S)$, where $ev$ is the natural evaluation map. It follows that as an element in $\hhh_K(Q \otimes \Omega, Q \otimes \Omega)$, $\theta_Q$ is given by $\Delta \circ ev$. Let $\{e_i\} , 1 \leq i \leq r$ be a basis for $Q$. Let $\{f_i\}$ be the basis of $Q^*$ dual to $\{e_i\}$. Let $\{u_i\}, 1 \leq i \leq n-r$ be a basis for $S$, and $\{v_i\}$ the basis for $S^*$ dual to $\{u_i\}$. What we have now discussed tells us the following lemma. \\

\begin{lemma}
With the notation just fixed, as an element of $\enn_K(Q \otimes \Omega) \bumpeq \enn(Q) \otimes \enn(\Omega) \bumpeq Q^* \otimes Q \bigotimes Q \otimes S^* \bigotimes Q^* \otimes S $ , $\theta_Q = \sum_{l_1,m_1,r_1} f_{m_1} \otimes e_{l_1} \bigotimes e_{m_1} \otimes v_{r_1} \bigotimes f_{l_1} \otimes u_{r_1}$ ($l_1$, $m_1$ running from $1$ to $r$, $r_1$ running from $1$ to $n-r$).
\end{lemma}

\begin{proof}
  $ev(e_i \otimes f_j \otimes u_k) = \delta_{ij}u_k$ and $\Delta(u_k) = \sum_{l=1}^r e_l \otimes f_l \otimes u_k$. Therefore,
$\theta_Q(e_i \otimes f_j \otimes u_k) =  \delta_{ij} \sum_{l=1}^r e_l \otimes f_l \otimes u_k$. On the other hand, $  f_{m_1} \otimes e_{l_1} \bigotimes e_{m_1} \otimes v_{r_1} \bigotimes f_{l_1} \otimes u_{r_1}(e_i \otimes f_j \otimes u_k) = \delta{im_1} \delta_{jm_1} \delta{kr_1} e_{l_1} \otimes f_{l_1} \otimes u_{r_1}$. This is nonzero iff $i=j=m_1$ and $k=r_1$. This proves the desired result.\\
\end{proof}

\subsection{Computing $\tilde{\ttt_k}(Q)$ for $k >1$}

This is done inductively. However, before we proceed, we need the following observations: \\

Observation 1. Let $V$ be any vector space over $K$. The following diagram (where arrows are $\sym^*V$-module homomorphisms) commutes:

\[
 \begin{CD}
     {V}^{\otimes d} \otimes \sym^*V @ > \alpha_i >> V^{\otimes d-1} \otimes \sym^*V\\
      @AAA                                               @AAA\\
     {\wedge}^dV  @> (-1)^{i-1}d >>              {\wedge}^{d-1}V
 \end{CD}
\] \\

Here, the vertical arrows are induced by the standard inclusions ${\wedge}^dV \subset V^{\otimes d}$ and ${\wedge}^{\otimes d-1}V \subset V^{d-1}$ respectively. By the standard inclusion ${\wedge}^kV \subset V^{\otimes k}$, we mean the map $ v_1 \wedge ... \wedge v_k \leadsto \sum_{\sigma \in S_k} \sn (\sigma) v_{\sigma(1)} \otimes ..... \otimes v_{\sigma(k)}$. $d$ denotes the Koszul differential in the bottom row. The map $\alpha_i$ in the top row is given by $(v_1 \otimes ... \otimes v_d) \bigotimes Y \leadsto \widehat{v_1 \otimes ...v_i...\otimes v_d} \bigotimes v_iY $. \\

Observation 2. If $x$ and $W$ are $P/N$ -representations (i.e, trivial $\sym^*\Omega$-modules) and if $\theta \in Hon_K(X,W)$ and $\tilde{\theta} \in \hhh_{\sym^* \Omega}(X \otimes \sym^*\Omega,W)$ is the element corresponding to $\theta$, then the following diagram commutes: \\

  \[
   \begin{CD}
   X \otimes \sym^* \Omega @> \tilde{\theta} >> W \\
    @V{(\theta \otimes \id)}VV             @VV{\id}V \\
   W \otimes \sym^* \Omega @> >>              W
 \end{CD}
 \] \\

To see what element in $\eee^k(Q, Q \otimes {\Omega}^{\otimes k})$ $\tilde{\ttt_k}(Q)$ is, we need to find vertical arrows so that the following diagram commutes, and describe the leftmost vertical arrow explicitly: \\

  \[
  \begin{CD}
   0 @> >> Q \otimes {\Omega}^{\otimes k} @> >> J_1(Q) \otimes {\Omega}^{\otimes k-1} @> >> .....@> >> J_1(Q)\\
 @.       @AAA                                  @AAA                                      @AAA       @AAA      \\
... @> >> Q \otimes {\wedge}^k \Omega \otimes \sym^* \Omega @> d >> Q \otimes {\wedge}^{k-1} \Omega \otimes \sym^* \Omega @> d >> ..... @ > d >> Q \otimes \sym^* \Omega
\end{CD}
\]  \\

 \[
  \begin{CD}
   J_1(Q) \otimes {\Omega}^{\otimes k-1} @> >> .....@> >> J_1(Q) @> >>Q @> >> 0\\
                              @AAA                                      @AAA       @AAA      @AA{\id}A \\
 Q \otimes {\wedge}^{k-1} \Omega \otimes \sym^* \Omega @> d >> ..... @ > d >> Q \otimes \sym^* \Omega @ > >> Q @ > >> 0
\end{CD}
\]  \\

Note that the bottom row in the above diagram is the Koszul complex. Observation 1 tells us that it suffices to find vertical arrows so that all squares in the following diagram commute and explicitly describe the leftmost nontrivial vertical arrow. \\

\[
 \begin{CD}
   0 @> >> Q \otimes {\Omega}^{\otimes k} @> >> J_1(Q) \otimes {\Omega}^{\otimes k-1} @> >> .....@> >> J_1(Q) \\
   @.       @AAA                                  @AAA                                      @AAA       @AAA      \\
  .. @> >> Q \otimes {\Omega}^{\otimes k} \otimes \sym^* \Omega @> \alpha_k >> Q \otimes {\Omega}^{\otimes k-1} \otimes \sym^* \Omega @> \alpha_{k-1} >> ..... @ >  >> Q \otimes \sym^* \Omega \\
    @.       @AAA                                  @AAA                                      @AAA       @AAA      \\
  .. @> >> Q \otimes {\wedge}^k \Omega \otimes \sym^* \Omega @> d >> Q \otimes {\wedge}^{k-1} \Omega \otimes \sym^* \Omega @> d >> ..... @ > d >> Q \otimes \sym^* \Omega
\end{CD}
 \] \\

\[
 \begin{CD}
    J_1(Q) \otimes {\Omega}^{\otimes k-1} @> >> .....@> >> J_1(Q) @> >>Q @> >> 0\\
                                    @AAA                                      @AAA       @AAA      @AA{\id}A \\
   Q \otimes {\Omega}^{\otimes k-1} \otimes \sym^* \Omega @> \alpha_{k-1} >> ..... @ >  >> Q \otimes \sym^* \Omega @ > >> Q @ > >> 0\\
                                    @AAA                                      @AAA       @AAA      @AA{\id}A \\
  Q \otimes {\wedge}^{k-1} \Omega \otimes \sym^* \Omega @> d >> ..... @ > d >> Q \otimes \sym^* \Omega @ > >> Q @ > >> 0
\end{CD}
 \] \\

Here, the bottom row is the Koszul complex and the map $\alpha_k:Q \otimes {\Omega}^{\otimes k} \otimes \sym^* \Omega \rar Q \otimes {\Omega}^{\otimes k-1} \otimes \sym^* \Omega$ is given by $ q \bigotimes w_1 \otimes ... \otimes w_k \bigotimes Y \leadsto q \bigotimes w_2 \otimes ... \otimes w_k \bigotimes w_1Y \text{ } \forall k$. Using this, we calculate $\tilde{\ttt_k}(Q)$ from $\tilde{\ttt_{k-1}}(Q)$ as follows: \\
$\tilde{\ttt_k}(Q) = (\tilde{\ttt_1}(Q) \otimes {\id_{\Omega}}^{\otimes k-1}) \circ \tilde{\ttt_{k-1}}(Q)$.
Suppose we have found vertical arrows which make the diagram below commute and suppose we have explicitly identified the leftmost nontrivial vertical arrow among these.  \\

\[
 \begin{CD}
   0 @> >> Q \otimes {\Omega}^{\otimes k-1} @> >> J_1(Q) \otimes {\Omega}^{\otimes k-2} @> >> .....@> >> J_1(Q)\\
@.        @A{\tilde{\theta}}AA               @A{\alpha}'AA                    @AAA @AAA \\
   .. @> >> Q \otimes {\Omega}^{\otimes k-1} \otimes \sym^* \Omega @> \alpha_{k-1} >> Q \otimes {\Omega}^{\otimes k-2} \otimes \sym^* \Omega @> \alpha_{k-2} >> ..... @ >  >> Q \otimes \sym^* \Omega
\end{CD}
 \] \\

\[
 \begin{CD}
    J_1(Q) \otimes {\Omega}^{\otimes k-2} @> >> .....@> >> J_1(Q) @> >>Q @> >> 0\\
               @A{\alpha}'AA                    @AAA @AAA @AA{\id}A\\
    Q \otimes {\Omega}^{\otimes k-2} \otimes \sym^* \Omega @> \alpha_{k-2} >> ..... @ >  >> Q \otimes \sym^* \Omega @ > >> Q @ > >> 0
\end{CD}
 \] \\

As we have already calculated $\tilde{\ttt_1}(Q)$ explicitly, we can draw the following commuting diagram:

\[
 \begin{CD}
   0 @> >> Q \otimes {\Omega}^{\otimes k} @> >> J_1(Q) \otimes {\Omega}^{\otimes k-1} @> >>  Q \otimes {\Omega}^{\otimes k-1} @> >> 0\\
@.       @A{\beta}AA          @A{\varphi}AA     @AA{\id}A\\
 .. @> >> Q \otimes \Omega \otimes {\Omega}^{\otimes k-1} \otimes \sym^* \Omega @> \alpha_{k} >> Q \otimes {\Omega}^{\otimes k-1} \otimes \sym^* \Omega @> >> Q \otimes {\Omega}^{\otimes k-1} @ > >> 0
\end{CD}
\]  \\

Here, $\beta = (\Delta \circ ev) \otimes \id_{{\Omega}^{\otimes k-1}}$, where $\Delta \circ ev = \theta_Q : Q \otimes \Omega \otimes \sym^* \Omega \rar  Q \otimes \Omega $. The diagram below now evidently commutes : \\

\[
  \begin{CD}
        Q \otimes {\Omega}^{\otimes k} @> >> J_1(Q) \otimes {\Omega}^{\otimes k-1}\\
 @AAA       @AA{\beta}A\\
 Q \otimes \Omega \otimes {\Omega}^{\otimes k-1} \otimes \sym^* \Omega @> \alpha_{k} >> Q \otimes {\Omega}^{\otimes k-1} \otimes \sym^* \Omega\\
 @A{\id_{\Omega} \otimes \theta \otimes \id_{\sym^* \Omega}}AA   @AA{\theta \otimes \id_{\sym^* \Omega}}A\\
 Q \otimes \Omega \otimes {\Omega}^{\otimes k-1} \otimes \sym^* \Omega @> \alpha_{k} >> Q \otimes {\Omega}^{\otimes k-1} \otimes \sym^* \Omega
 \end{CD}
\] \\

Moreover, by Observation 2, the diagram below commutes: \\

\[
  \begin{CD}
     J_1(Q) \otimes {\Omega}^{\otimes k-1} @> >>  Q \otimes {\Omega}^{\otimes k-1} @> >> J_1(Q) \otimes {\Omega}^{\otimes k-2}\\
@A{\varphi}AA            @A{\id}AA    @AA{\id}A\\
Q \otimes {\Omega}^{\otimes k-1} \otimes \sym^* \Omega @> >>  Q \otimes {\Omega}^{\otimes k-1} @> >> J_1(Q) \otimes {\Omega}^{\otimes k-2}\\
@A{\theta \otimes \id_{\sym^* \Omega}}AA    @A{\theta}AA   @AA{\alpha}'A\\
Q \otimes {\Omega}^{\otimes k-1} \otimes \sym^* \Omega @> \id >> Q \otimes {\Omega}^{\otimes k-1} \otimes \sym^* \Omega @ > \alpha_{k-1} >> Q \otimes {\Omega}^{\otimes k-2} \otimes \sym^* \Omega
\end{CD}
\]  \\

From the above two diagrams, we see that the diagram below, diving us $\tilde{\ttt_k}(Q)$ commutes: \\

\[
  \begin{CD}
   0 @> >> Q \otimes {\Omega}^{\otimes k} @> >> J_1(Q) \otimes {\Omega}^{\otimes k-1} @> >>J_1(Q) \otimes {\Omega}^{\otimes k-2} @> >> .....\\
@. @A{\beta \circ (\id_{\Omega} \otimes \theta \otimes \id_{\sym^* \Omega})}AA @A{\varphi \circ (\theta \otimes \id_{\sym^* \Omega})}AA @A{\alpha}'AA @AAA \\  .. @> >> Q \otimes {\Omega}^{\otimes k} \otimes \sym^* \Omega @> \alpha_{k} >> Q \otimes {\Omega}^{\otimes k-1} \otimes \sym^* \Omega @> \alpha_{k-1} >> Q \otimes {\Omega}^{\otimes k-2} \otimes \sym^* \Omega @> \alpha_{k-2} >>....
\end{CD}
\]\\

\[
  \begin{CD}
   J_1(Q) \otimes {\Omega}^{\otimes k-1} @> >>J_1(Q) \otimes {\Omega}^{\otimes k-2} @> >> .....@> >> Q @> >> 0\\
 @A{\varphi \circ (\theta \otimes \id_{\sym^* \Omega})}AA @A{\alpha}'AA @AAA @A{\id}AA\\  Q \otimes {\Omega}^{\otimes k-1} \otimes \sym^* \Omega @> \alpha_{k-1} >> Q \otimes {\Omega}^{\otimes k-2} \otimes \sym^* \Omega @> \alpha_{k-2} >>.... @ >  >>  Q @ > >> 0
\end{CD}
\]\\

This tells us what $\tilde{\ttt_k}(Q)$ is as an element of $\hhh_K( Q \otimes {\wedge}^k \Omega, Q \otimes {\Omega}^{\otimes k})$. Since we have the standard inclusion $i: {\wedge}^k \Omega \rar {\Omega}^{\otimes k}$ which gives rise to $i: Q \otimes  {\wedge}^k \Omega \rar Q \otimes {\Omega}^{\otimes k}$, we can write $\tilde{\ttt_k}(Q) = \gamma_k \circ i$ , where $\gamma_k \in \hhh_K( Q \otimes {\Omega}^{\otimes k}, Q \otimes {\Omega}^{\otimes k})$. The preceding discussion tells us that $ \gamma_k = (( \Delta \circ ev) \otimes \id_{{\Omega}^{\otimes k-1}}) \circ \id_{\Omega} \otimes \gamma_{k-1}$. From this, we obtain the next lemma (notation as in the previous subsection).

\begin{lemma}
Identifying $\enn_K( Q \otimes {\Omega}^{\otimes k})$ with $\enn_K(Q) \bigotimes {{\Omega}^*}^{\otimes k} \bigotimes {\Omega}^{\otimes k}$ , we have : $$ \gamma_k = \sum_{l_1,...,l_k;m_1,....,m_k;r_1,...,r_k} (f_{m_1} \otimes e_{l_1}) \circ ... \circ (f_{m_k} \otimes e_{l_k}) \bigotimes (e_{m_1} \otimes v_{r_1}) \otimes .....\otimes (e_{m_k} \otimes v_{r_k})$$ $$ \bigotimes (f_{l_1} \otimes u_{r_1}) \otimes.....\otimes (f_{l_k} \otimes u_{r_k})$$ . Here, the $l_i , 1 \leq i \leq k $ and the $m_i , 1 \leq i \leq k$ run from $1$ to $r$ , while the $r_i , 1 \leq i \leq k$ run from $1$ to $n-r$.
\end{lemma}

\begin{proof} By induction on $k$. The base case is the formula for $\theta_Q$ and the induction step is the formula $ \gamma_k = (( \Delta \circ ev) \otimes \id_{{\Omega}^{\otimes k-1}}) \circ \id_{\Omega} \otimes \gamma_{k-1}$. \\
\end{proof}

Having computed $\tilde{\ttt_k}(Q)$ we compute $\ttt_k(Q)$ . For this, we need to show that $\ttt_k(Q) = (tr \otimes \id)_* \tilde{\ttt_k}(Q)$ where $\tilde{\ttt_k}(Q) \in \enn(Q) \otimes \hhh_K({\wedge}^k \Omega, {\Omega}^{\otimes k})$ and $tr: \enn(Q) \rar K$ is the trace map. This we shall prove in the appendix. Calculating $\ttt_k(Q)$ is then easy : in the formula in the previous lemma, we see that $ (f_{m_1} \otimes e_{l_1}) \circ ... \circ (f_{m_k} \otimes e_{l_k})(e_i) = \delta_{im_k} \delta_{l_km_{k-1}}....\delta_{l_2m_1}e_{l_1}$. From this, we see that $(f_{m_1} \otimes e_{l_1}) \circ ... \circ (f_{m_k} \otimes e_{l_k})$ has trace $1$ iff $ m_k =l_1 , l_k=m_{k-1} , .....,l_2=m_1$ and has trace $0$ otherwise. From this it follows that if $i:{\wedge}^k \Omega \rar {\Omega}^{\otimes k}$ is the natural inclusion, $\ttt_k(Q)$ is given by $ \mu_k \circ i$ where $\mu_k \in \hhh_K({\Omega}^{\otimes k}, {\Omega}^{\otimes k})$ is as described in the following lemma  \\

\begin{lemma}
  Identifying $\enn_K({\Omega}^{\otimes k})$ with ${{\Omega}^*}^{\otimes k} \bigotimes {\Omega}^{\otimes k}$ we have
$$ \mu_k = \sum_{l_1,....,l_k;r_1,...,r_k} (e_{l_2} \otimes v_{r_1}) \otimes ....\otimes (e_{l_k} \otimes v_{r_{k-1}}) \otimes (e_{l_1} \otimes v_{r_k})\bigotimes (f_{l_1} \otimes u_{r_1})\otimes.... \otimes (f_{l_k} \otimes u_{r_k})$$
$$ = \sum_{m_1,...,m_k;r_1,...,r_k}(e_{m_1} \otimes v_{r_1}) \otimes ....\otimes (e_{m_k} \otimes v_{r_{k}}) \bigotimes (f_{m_k} \otimes u_{r_1})\otimes (f_{m_1} \otimes u_{r_2})\otimes.... \otimes (f_{m_{k-1}} \otimes u_{r_k})$$ \\
\end{lemma}

As a consequence, the basis element $f_{i_1} \otimes .... \otimes f_{i_k} \bigotimes u_{j-1} \otimes... \otimes u_{j_k}$ of ${\Omega}^{\otimes k}$ is mapped by $\ttt_k(Q)$ to $f_{i_k} \otimes f_{i_1} \otimes .... \otimes f_{i_{k-1}} \bigotimes u_{j-1} \otimes... \otimes u_{j_k}$ Therefore, if we identify $\enn_K({\Omega}^{\otimes k})$ with ${Q^*}^{\otimes k} \otimes S^{\otimes k}$, $\ttt_k(Q)$ ca be thought of as $(k \text{ } k-1 \text{ } k-2 \text{ }.. \text{ } 2 \text{ } 1) \otimes \id_{ S^{\otimes k}}$ where $(k \text{ } k-1 \text{ } k-2 \text{ }.. \text{ } 2 \text{ } 1)$ is the $k$-cycle acting on ${Q^*}^{\otimes k}$ by the usual action of $S_k$ on $V^{\otimes k}$ for a vector space $V$. \\

\section{ Proofs of Theorems 1 and 2}

First, we must explicitly describe what $\sigma_* \ttt_k(Q)$ looks like as an element of $\hhh_K({\wedge}^k \Omega, {\Omega}^{\otimes k})$ where $\sigma \in S_k$ acts by $v_1 \otimes.... \otimes v_k \leadsto v_{\sigma(1)} \otimes.... \otimes v_{\sigma(k)}$ . That ${\sigma}_*\ttt_k(Q) = \sigma \circ \ttt_k(Q)$ , where on the right hand side $\sigma$ is the endomorphism of ${\Omega}^{\otimes k}$ just described is seen from the commutative diagram below where the middle row describes $\ttt_k(Q)$ as an element of $\eee^k(K, {\Omega}^{\otimes k})$ ,and where the bottom row is a Koszul complex. \\

\[
 \begin{CD}
  0@> >>{\Omega}^{\otimes k} @> >>Z_1@> >>...@> >>Z_k@> >>K @> >>0\\
 @.  @A{\sigma}AA    @AAA @AAA @AAA @AA{\id}A\\
  0@> >>{\Omega}^{\otimes k} @> >>Y_1@> >>...@> >>Y_k@> >>K @> >>0\\
@. @AAA  @AAA  @AAA  @AAA  @AA{\id}A\\
...@> >>{\wedge}^k \Omega  \otimes \sym^* \Omega @> >>{\wedge}^{k-1} \Omega  \otimes \sym^* \Omega @> >>.....@> >>\sym^* \Omega @> >>K @> >>0 \\
\end{CD}
\]

Note that $\ttt_k(Q) = \mu_k \circ i$ by Lemma 10 where $i:{\wedge}^k \Omega \rar {\Omega}^{\otimes k}$ is the natural inclusion. Therefore, $\sigma_*\ttt_k(Q) = \sigma \circ \mu_k \circ i$. However, we have a projection $p:{\Omega}^{\otimes k} \rar {\wedge}^k \Omega$ given by $v_1 \otimes ... \otimes v_k \leadsto \sum_{\sigma \in S_k} \frac{1}{k!} \sn(\sigma)(v_{\sigma(1)} \otimes ... \otimes v_{\sigma(k)})$. Note that $ \alpha \circ i = 0$ iff $\alpha \circ p =0$, for any $\alpha \in \enn_K({\Omega}^{\otimes k})$. Now observe that if $\sigma: {\Omega}^{\otimes k} \rar {\Omega}^{\otimes k}$ is given by  $\sigma \otimes \sigma:\enn({Q^*}^{\otimes k}) \otimes \enn(S^{\otimes k})$ when we identify $\enn_K({\Omega}^{\otimes k})$ with $\enn_K({Q^*}^{\otimes k}) \otimes \enn_K(S^{\otimes k})$. It follows that as an element of  $\enn_K({Q^*}^{\otimes k}) \otimes \enn_K(S^{\otimes k})$, $\ttt_k(Q)$ is given by $\frac{1}{k!} (\tau_k \otimes \id_{S^{\otimes k}}) \circ \sum_{\omega \in S_k} \sn(\omega)( \omega \otimes \omega)$ as an element of the image of $\hhh_K({\wedge}^k \Omega , {\Omega}^{\otimes k})$ in $\enn_K( {\Omega}^{\otimes k})$. As the action of $S_k$ on $V^{\otimes k}$ given by  $v_1 \otimes ... \otimes v_k \leadsto v_{\sigma(1)} \otimes ... \otimes v_{\sigma(k)}$ is a right action, we see that $\ttt_k(Q) = \sum_{\omega \in S_k} \sn(\omega)( \omega \tau_k \otimes \omega)$. Therefore, $\sigma_*\ttt_k(Q) =  \frac{1}{k!} \sum_{\omega \in S_k} \sn(\omega)( \omega \tau_k \otimes \omega)( \sigma \otimes \sigma) = \frac{1}{k!} \sum_{\omega \in S_k} \sn(\omega)( \omega \tau_k \sigma \otimes \omega \sigma) = \frac{1}{k!} \sum_{\omega \in S_k} \sn(\omega)( \omega \sigma \sigma^{-1} \tau_k \sigma \otimes \omega \sigma) = \frac{1}{k!} \sum_{\omega \sigma \in S_k} \sn(\omega \sigma) \sn(\sigma)( \omega \sigma \sigma^{-1} \tau_k \sigma \otimes \omega \sigma) =  \frac{1}{k!} \sum_{\beta \in S_k} \sn(\beta) \sn(\sigma)( \beta \sigma^{-1} \tau_k \sigma \otimes \beta) = \frac{1}{k!} \sum_{\beta \in S_k} \sn(\beta)( \beta \otimes \beta) \sigma^{-1} \tau_k \sigma \otimes  \id_{S^{\otimes k}})$ . (Here $\tau_k$ denotes the $k$-cycle $(k \text{ } k-1 \text{ }....\text{ } 2 \text{ } 1)$.) This therefore, proves the following lemma: \\

\begin{lemma}
$ \sigma_*\ttt_k(Q) = \sn(\sigma) \sigma^{-1} \tau_k \sigma \otimes  \id_{S^{\otimes k}} \circ i $ where $i:{\wedge}^k \Omega \rar {\Omega}^{\otimes k}$ is the natural inclusion. \\
\end{lemma}

Therefore, if $\sum a_{\sigma}\sn(\sigma) \sigma \in K(S_k)$ , we have $\sum a_{\sigma}\sn(\sigma) \sigma_* (\ttt_k(Q)) = [(\sum a_{\sigma}\sn(\sigma) \sigma^{-1} \tau_k \sigma) \otimes \id_{S^{\otimes k}}] \circ i $. It follows from the corollary to Lemma  describing the projector giving $\ttt_k(\alpha_l(V))$ from $\ttt_k(V)$ that \\

\begin{lemma}
$\ttt_k(\alpha_l(Q)) = \sum_{j=1}^n  \sum_{\text{$\sigma$ has $j-1$ descents}} ({a^{l,j}}_k \sigma \tau_k \sigma^{-1}) \otimes \id_{S^{\otimes k}}$
\end{lemma}

After all $\sum_{j=1}^n  \sum_{\text{$\sigma$ has $j-1$ descents}} \sn(\sigma) a_{\sigma} \sigma^{-1}$ is the operator ${e^{(l)}}_k$ for the graded commutative Hopf-algebra $T^*V$. In fact, $\sum_{j=1}^n  \sum_{\text{$\sigma$ has $j-1$ descents}} a_{\sigma} \sigma$ is the operator ${e^{(l)}}_k$ for the co-commutative ordinary Hopf-algebra $TV$. We henceforth denote this idempotent by  ${\tilde{e}^{(l)}}_k$. Let $*$ denote the conjugation action of $KS_k$ on itself, i.e, if $a \in S_k$ and $b \in KS_k$ then $a*b = aba^{-1}$, and $\sum c_gg*h = \sum c_gghg^{-1} , h \in KS_k$. Then, we have $\mu_*\ttt_k(Q) = \mu * \tau_k \otimes \id_{S^{\otimes k}} \forall \mu \in KS_k $  and therefore that  $\ttt_k(\alpha_l(Q)) = {\tilde{e}^{(l)}}_k * \tau_k \otimes \id_{S^{\otimes k}}$. Note that $*$ is a left action. \\
Also observe that $(\enn_K(\Omega^{\otimes k}))^{P/N} = (\enn_K(Q^*) \otimes \enn_K(S))^{Gl(Q) \times Gl(S)} $. It is a theorem of Weyl that $KS_k \otimes KS_k \rar (\enn_K(Q^*) \otimes \enn_K(S))^{Gl(Q) \times Gl(S)}$ is surjective. We therefore, need to find the elements $ \alpha \in KS_k$ such that $\alpha \otimes \id \circ p =0 $ as an element of $ \enn_K(Q^*) \otimes \enn_K(S))^{Gl(Q) \times Gl(S)}$ . The answer depends on the ranks of $Q$ and $S$. Before we proceed we need a small digression. \\

\subsection{A lemma and some corollaries}

\begin{lemma}
Let $G$ be a finite group, and let $\chi: G \rar {\mathbb C}^*$ be a $1$-dimensional representation of $G$. Then, if  $\beta \in {\mathbb C}(G)$, $\sum_{g \in G}\chi(g) (g \otimes g) (\beta \otimes \id) =0$ in ${\mathbb C}(G \times G) ={\mathbb C}(G) \otimes {\mathbb C}(G)$ iff $\beta =0$.
\end{lemma}

\begin{proof} If $\beta =0$ then clearly $\sum_{g \in G}\chi(g) (g \otimes g) (\beta \otimes \id) =0$. For the implication in the opposite direction, let us see what $\sum_{g \in G}\chi(g) (g \otimes g)$ does to ${\mathbb C}(G \times G) = \oplus \enn(V_x \otimes V_y)$ where the $V_x$ are the irreducible representations of $G$. Let $e_i$ be a basis for $V_x$ and let $f_j$ be a basis of $V_y$. Consider a matrix representation of ${\mathbb C}(G)$ associated to this choice of basis. Then, $\sum_g \chi(g) (g \otimes g) (e_i \otimes f_j) = \sum_g \sum_{k,l} {g_{ik}}^x {g_{jl}}^y \chi(g) (e_k \otimes f_l) = \sum_{k,l} (e_k \otimes f_l) (\sum_g \chi(g){g_{ik}}^x {g_{jl}}^y = \sum_{k,l} (e_k \otimes f_l) (\sum_g {g_{ik}}^z {g_{jl}}^y)$ , where $V_z = V_x \otimes \chi$. \\
Note that $\sum_g (g \otimes g) \in  \enn(V_z \otimes V_{y})$ is a $G$-module homomorphism. In fact, $G$ acts trivially on $(\sum_g g \otimes g ).(V_z \otimes V_{y})$. Thus, $\sum_g (g \otimes g)$ acts as a projection to the trival part of $V_z \otimes V_{y})$. Note that $V_z \otimes V_{y}$ has a contains precisely $\langle \chi_z , \bar{\chi}_{y}  \rangle $ copies of the trivial representation of $G$. In particular, it contains one copy of the trivial representation of $G$ iff $V_z$ and $V_{y}$ are dual representations. In that case, the projection to that copy of the trivial representation is given by $v \otimes w \leadsto w(v) \sum e_i \otimes f_i$. Here, $\{e_i\}$ is a basis for $V_z$ and $\{f_i\}$ is the basis for $V_{y}$ dual to $\{e_i\}$. This tells us that $\sum_g {g_{ik}}^z {g_{jl}}^y = \delta_{y \bar{z}} \delta_{ij} \delta_{kl}$. \\
Therefore, in $\enn(V_x \otimes V_{y})$ , if $V_z$ is not dual to $V_{y}$, then  $\sum_{g \in G}\chi(g) (g \otimes g) = 0$. If $V_z$ is dual to $V_{y}$, then if  $\{e_i\}$ is a basis for $V_z$ and $\{f_i\}$ is the basis of $V_{y}$ dual to $\{e_i\}$, and if $\{\tilde{e}_i\}$ is the basis of $V_x$ corresponding to $\{e_i\}$, then with respect to the ordered basis $\tilde{e}_1 \otimes f_1 , \tilde{e}_2 \otimes f_1,....,\tilde{e}_d \otimes f_1,\tilde{e}_1 \otimes f_2,......,\tilde{e}_d \otimes f_2, .......,\tilde{e}_1 \otimes f_d, .....,\tilde{e}_d \otimes f_d$ of $V_x \otimes V_{y}$, $\sum_{g \in G}\chi(g) (g \otimes g)$ corresponds to the matrix $M$ such that $M_{ij} = 1 \text{ if } i ,j \in \{ kd +k+1 | 0 \leq k \leq d-1\}$. $ M_{ij} =0 $ otherwise. On the other hand , $\beta \otimes \id$ in $\enn(V_x \otimes V_{y})$ is given by a block diagonal matrix each of whose diagonal blocks is the matrix representing $\beta$ in $\enn(V_x)$. This proves the desired lemma. \\

\end{proof}

In fact, in the above proof, we have done a bit more:\\

\begin{lemma}
Let $G$ be a finite group, and let $\chi: G \rar {\mathbb C}^*$ be a $1$-dimensional representation of $G$. Let $V_x$ and $V_{y}$ be irreducible representations of $G$ such that $V_x \otimes \chi$ is dual to $V_{y}$. Then, if  $\beta \in {\mathbb C}(G)$, $\sum_{g \in G}\chi(g) (g \otimes g) (\beta \otimes \id) =0$ in $\enn(V_x \otimes V_{y})$ iff $\beta =0$ in $\enn(V_x)$. \\
\end{lemma}

In our problem , the group in question is $S_k$. We note that the irreducible representations of $S_k$ over $\mathbb C$ can be realised over $\mathbb Q$ and hence over any field of characteristic $0$. The previous two lemmas thus hold for $KS_k$ where $K$ is a field of characteristic $0$. We also recall that the irreducible representations of $S_k$ are indexed by partitions $\lambda$ of $k$. They are self-dual, and $V_{\lambda} \otimes Alt = V_{\bar{\lambda}}$ , where $\bar{\lambda}$ is the partition conjugate to ${\lambda}$. \\
These lemmas give us the following framework in which we can view our problem: Let ${\mathbb S}_{\lambda}$ denote the Schur-functor associated with the partition $\lambda$ of $k$ i.e, if $V$ is any vector space ${\mathbb S}_{\lambda}(V) = V^{\otimes k} \otimes_{KS_k} V_{\lambda}$ where $V_{\lambda}$ is the irreducible representation of $S_k$ corresponding to the partition $\lambda$. We know that if $V$ is a vector space of rank $m$ ,  ${\mathbb S}_{\lambda}(V) =0$ iff $\lambda$ has more than $m$ parts. Therefore if $Q$ has rank $r$, then ${\mathbb S}_{\lambda}(Q) = 0$ iff $|\lambda| > r$ and  ${\mathbb S}_{\bar{\lambda}}(S)=0$ iff $|\bar{\lambda}| > n-r$. Moreover, if $\lambda$ and $\mu$ are two partitions of $k$ , then $V^{\otimes k} \otimes W^{\otimes k} \otimes_{K(S_k \times S_k)} V_{\lambda} \otimes V_{\mu} =  {\mathbb S}_{\lambda}(V) \otimes  {\mathbb S}_{\mu}(W)$. If $\gamma \in K(S_k \times S_k) \neq 0 $ in $\enn(V_{\lambda} \otimes V_{\mu})$, then $K(S_k \times S_k). \gamma $ contains $V_{\lambda} \otimes V_{\mu}$. Therefore, $V^{\otimes k} \otimes W^{\otimes k} \otimes_{K(S_k \times S_k)} \gamma$ contains ${\mathbb S}_{\lambda}(V) \otimes {\mathbb S}_{\mu}(W)$. Lemma  therefore , says the following: \\

\begin{lemma}
If the rank of $Q$ is $r$ and that of $S$ is $n-r$ , then $\sum_{\sigma} \sn(\sigma) ( \sigma \otimes \sigma) (\beta \otimes \id) = 0$ as an element of $\hhh_K( {\Omega}^{\otimes k} , {\Omega}^{\otimes k})$ iff $\beta = 0$ as an element of $\enn(V_{\lambda})$ for all partitions $\lambda$ such that $|\lambda| \leq r$ and $| \bar{\lambda}| \leq n-r$
\end{lemma}

\begin{proof}
Let $ \gamma =  \sum_{\sigma} \sn(\sigma) ( \sigma \otimes \sigma) (\beta \otimes \id)$. Then, $\gamma = 0$ in $\enn(V_{\lambda} \otimes V_{\mu})$ if $\mu \neq \bar{\lambda}$. Therefore, $\gamma$ kills ${\mathbb S}_{\lambda}(Q^*) \otimes {\mathbb S}_{\mu}(S)$ whenever  $\mu \neq \bar{\lambda}$. On the other hand, if $\gamma \neq 0$ in $\enn(V_{\lambda} \otimes V_{\bar{\lambda}})$,  then, ${\Omega}^{\otimes k} . \gamma$ contains a copy of ${\mathbb S}_{\lambda}(Q^*) \otimes {\mathbb S}_{\bar{\lambda}}(S)$. The desired lemma follows immediately. \\
\end{proof}

Thus, if $Q$ is of rank $r$, the annihilator in $KS_k$ of $\ttt_k(\alpha_l(Q))$ is precisely the space of all $\sum_g c_g g$ so that $(\sum_g c_g g^{-1})* {\tilde{e}^{(l)}}_k * \tau_k = 0 $ as an element of $\enn(V_{\lambda})$ for all $\lambda$ satisfying $|\lambda| \leq r$ and $| \bar{\lambda}| \leq n-r$. For the time being, let us assume that $n$ is large enough , so that we do not have to bother about the condition  $| \bar{\lambda}| \leq n-r$. Then , $I(k,r,l) = \{ \sum_g c_g g | (\sum_g c_g g^{-1})* {\tilde{e}^{(l)}}_k * \tau_k = 0 \text{ as an element of } \enn(V_{\lambda}) \text{ } \forall |\lambda| \leq r\}$ .On the other hand, $I(k, r-1, l) = \{ \sum_g c_g g | (\sum_g c_g g^{-1})* {\tilde{e}^{(l)}}_k * \tau_k = 0 \text{ as an element of } \enn(V_{\lambda}) \text{ } \forall |\lambda| \leq r-1\}$. It is now clear that $I(k,r,l) \subseteq I(k, r-1,l)$ . What we need to show is that for  fixed $l$ and $r$, there exists a $k$ for which the above inclusion is strict. Let us denote by $P_r$ the projection from $KS_k$ to $\oplus_{|\lambda| \leq r} \enn(V_{\lambda}) \text{ } 1 \leq r \leq k$. Then $I(k,r,l) = \{  \sum_g c_g g | (\sum_g c_g g^{-1})* P_r({\tilde{e}^{(l)}}_k * \tau_k) =0\}$. and $\dim(I(k,r,l)) = \dim( \langle {\tilde{e}^{(l)}}_k * \tau_k \rangle ) - \dim( \langle  P_r({\tilde{e}^{(l)}}_k * \tau_k)  \rangle )$ (if $k < n-r$) where for a given $\alpha \in KS_k$, $ \langle  \alpha \rangle $ denotes the subspace of $KS_k$ spanned by $ \beta * \alpha  , \beta \in S_k$. We must thus show that fixed $l$ and $r$ , there exists a $k$ so that $\dim( \langle  P_r({\tilde{e}^{(l)}}_k * \tau_k)  \rangle ) > \dim( \langle  P_{r-1}({\tilde{e}^{(l)}}_k * \tau_k)  \rangle )$. We will show that for a fixed $l$, there are infinitely many $r$ for which such a $k$ exists. This we will do using a simple counting argument. But we need some further lemmas before we can proceed. \\

\begin{lemma}

${\tilde{e}^{(l)}}_k * \tau_k = {\tilde{e}^{(l-1)}}_{k-1} * \tau_k$ where $S_{k-1} \subset S_k$ is embedded as the subgroup fixing $k$.

\end{lemma}

\begin{proof}

Let $\alpha$ be a permutation of $\{1,2,3,...,k-1\}$ with $j-1$ descents. Then, among the permutations $\alpha, \alpha \tau_k ,...., \alpha {\tau_k}^{k-1}$ , we see that $j$ of the permutations have $j-1$ descents, while the remaining $k-j$ have $j$ descents. For, $\alpha {\tau_k}^i$ has $j$ descents or $j-1$ descents depending on whether $\alpha(k-i) < \alpha(k-i+1)$ or not, for $2 \leq i \leq k-1$. For $j-1$ such $i$ ,  $\alpha(k-i) > \alpha(k-i+1)$ (corresponding to the descents of $\alpha$). These $j-1$ elements together with $\alpha$ have $j-1$ descents. The remaining $k-j$ permutations have $j$ descents. As $ {\tau_k}^i \tau_k {\tau_k}^{-i} = \tau_k$, the coefficient of $\alpha \tau_k \alpha^{-1}$ in $ {\tilde{e}^{(l)}}_k * \tau_k$ is given by $j{a_k}^{l,j} + (k-j){a_k}^{l,j+1}$, since among the elements  $\alpha, \alpha \tau_k ,...., \alpha {\tau_k}^{k-1}$ , those with $j-1$ descents contribute $ {a_k}^{l,j}$ and those with $j$ descents contribute ${a_k}^{l,j+1}$ to the coefficient of  $\alpha \tau_k \alpha^{-1}$ in $ {\tilde{e}^{(l)}}_k * \tau_k$. The desired lemma follows from observing that   $j{a_k}^{l,j} + (k-j){a_k}^{l,j+1} = j{a_{k-1}}^{l-1,j} $, since $j \binom{X-j+k}{k} +(k-j) \binom{X-j-1-k}{k} = X \binom{X-j-1-k}{k-1}$. \\

\end{proof}

Also observe that the stabiliser of ${\tau_k}$ under conjugation is the cyclic usubgroup generated by $\tau_k$. Thus, $S_{k-1}$ acts freely on the conjugates of $\tau_k$ and $\beta * \tau_k =0$ for some $\beta \in KS_{k-1}$ iff $\beta =0$. It follows from this remark and the above lemma that $\dim( \langle {\tilde{e}^{(l)}}_k * \tau_k \rangle )$ is the dimension of the representation $KS_{k-1}.{\tilde{e}^{(l-1)}}_{k-1}$ of $KS_{k-1}$. It is a result in Loday[2] that this space has dimension equal to the coefficient of $q^{l-1}$ in $q(q+1)...(q+k-2)$. \\
On the other hand, look at $\dim( \oplus_{|\lambda| \leq r} \enn(V_{\lambda})$ for a fixed $r$. Note that if $\lambda: k= \lambda_1+... + \lambda_{r'}$ is a partition of $k$, and if $\Pi$ denotes the product of the hook lengths of the Young diagram corresponding to ${\lambda}$, then $\dim(V_{\lambda}) = \frac{k!}{\Pi} \leq \frac{k!}{\lambda_1!\lambda_2!...\lambda_{r'}!}$. Thus, $\dim(\enn(V_{\lambda})) \leq {\frac{k!}{\lambda_1!\lambda_2!...\lambda_{r'}!}}^2$. Hence, $\dim(( \oplus_{|\lambda| \leq r} \enn(V_{\lambda}) \leq \sum_{\lambda_1 + ...+ \lambda_r = k ; \lambda_i \geq 0} {\frac{k!}{\lambda_1!\lambda_2!...\lambda_{r}!}}^2 \leq (\sum_{\lambda_1 + ...+ \lambda_r = k ; \lambda_i \geq 0} \frac{k!}{\lambda_1!\lambda_2!...\lambda_{r}!})^2 = r^{2k}$. Therefore, for a fixed $r$, $\dim( \langle  P_r({\tilde{e}^{(l)}}_k * \tau_k)  \rangle ) \leq \dim(( \oplus_{|\lambda| \leq r} \enn(V_{\lambda}) \leq r^{2k}$.  On the other hand $\dim( \langle {\tilde{e}^{(l)}}_k * \tau_k \rangle ) = \text{ coefficient of } q^{l-1} \text{ in } q(q+1)...(q+k-2) \geq  \frac{(k-2)!}{(l-2)!}$. \\
Now, for sufficiently large $k$, $\frac{(k-2)!}{(l-2)!} > r^{2k}$. To see this we need to find $k$ large enough so that $\ln( (k-2)!) - \ln( (l-2)!) > 2k \ln r$. $\ln( (k-2)!) > (k-2)\ln( k-2) -(k-3)$. We therefore, only need to find $k$ large enough so that $(k-2)\ln(k-2) > k-3 + \ln((l-2)!) + (k-2)\ln( r^2) + 2\ln(r^2)$. Put $D = \ln( r^4 (l-2)!)$. We then need $k$ so that $ (k-2)\ln(k-2) > k-3 + D + (k-2)\ln( r^2)$. Certainly, $\exists N \in {\mathbb N}$ so that $N(k-2) > (k-3) +D$ (if we pick $N > D+1$ and $k >3$ for instance). If $k-2 > e^{N}r^2$, then we see that $ (k-2)\ln(k-2) > k-3 + D + (k-2)\ln( r^2)$. Certainly, $k > e^{N+1}r^2$ would do for our purposes. Thus, if $l$ and $r$ are fixed, we have shown that there is a constant $C$ so that when $k > Cr^2$, then  $\dim( \langle {\tilde{e}^{(l)}}_k * \tau_k \rangle ) > \dim( \langle  P_r({\tilde{e}^{(l)}}_k * \tau_k)  \rangle )$. If $l=2$, in particular, we need $(k-2)\ln(k-2) > k-3 + (k-2)\ln( r^2) + 2\ln(r^2)$. We see that this happens if $k-2 > 7r^2$. \\
Note that $\dim( \langle {\tilde{e}^{(l)}}_k * \tau_k \rangle ) > \dim( \langle  P_r({\tilde{e}^{(l)}}_k * \tau_k)  \rangle )$ tells us that $ \exists s \geq r \text{ so that } \dim( \langle  P_s({\tilde{e}^{(l)}}_k * \tau_k)  \rangle ) < \dim( \langle  P_{s+1}({\tilde{e}^{(l)}}_k * \tau_k)  \rangle )$. Thus , $\exists s \geq r$ so that $I(k, s+1,l) \subsetneq I(k,s,l)$. This infact completes the proof of Theorem 1 and the bulk of the proof of Theorem 2. To complete the proof of Theorem 2, we make some observations: \\

 Observation 1: By Lemma 16 $\tau_k = \sum_{l \geq 2} {\tilde{e}^{(l-1)}}_{k-1} * \tau_k = \sum_{l \geq 2} {\tilde{e}^{(l)}}_k * \tau_k \implies \ttt_k(Q) = \sum_{l \geq 2} \ttt_k(\alpha_l(Q)) \implies \ttt_k(\alpha_1(Q)) =0 \text{ } \forall k \geq 2$. \\

Observation 2: Since $\oplus \ttt_k: K(X) \otimes {\mathbb Q} \rar \oplus \Hm^k(X , {\Omega}^{\otimes k})$ is a ring homomorphism, is follows that $\ttt_k( {\alpha_1(Q)}^2) = 0$ if $k \neq 2$. \\

If $f: G(s+1, N) \rar G(s,M)$ is a morphism, then one sees that $f^*(\alpha_2(Q')) = A{\alpha_1(Q)}^2 + B\alpha_2(Q)$. By the Observation 2, $\ttt_k(f^*(\alpha_2(Q'))) =  B\ttt_k(\alpha_2(Q))$. If $B \neq 0$, one sees that $I(k,s,2) \subseteq I(k,s+1,2)$ (a contradiction). This finally proves Theorem 2. \\

To prove Theorem 4, we need the following lemma from which Theorem 4 follows immediately.   \\

\begin{lemma}
X a smooth (projective) scheme. Suppose that $[V] \in K(X)$ is given by $[V] = \sum a_i[V_i]$, where $V_i$'s are of rank $\leq r$. Then, $I(k,r,l)$ annihilates $\ttt_k([V])$.
\end{lemma}

\begin{proof}
$ \exists N \in {\mathbb N} $ so that $ \forall M > N$ there exist surjections ${\mathbb G}_i \rar V_i(m)$ where ${\mathbb G}_i$ is a free $\calo_X$ module $\forall i$. This is equivalent to saying that there exist morphisms $f_i: X \rar G(rank(V_i) , K_i)$ so that $V_i(m) = {f_i}^*Q_i$, $Q_i$ the universal quotient bundle of  $G(rank(V_i) , K_i) \text{ } \forall i$. Thus, $I(k,r,l)$ kills $\ttt_k(\alpha_l(V_i \otimes \calo(m))) \forall m > N , \forall i$ . It suffices to show that $I(k,r,l)$ kills $\ttt_k(\alpha_l(V_i)) \forall i $. For this, we note that $\oplus \ttt_k(\calo(1)) = e^{\ttt_1(\alpha_1(\calo_1))}$, with the understanding that ${\ttt_1(\alpha_1(\calo_1))}^{D+1} =0$ where $D$ is the dimension of the ambient projective space. Thus, $\oplus \ttt_k(\calo(m)) = e^{m\ttt_1(\alpha_1(\calo_1))}$. Since the Vandermonde determinant $\Delta(N+1,..,N+D+1) \neq 0$, we can find a linear combination $W$ of $\calo(N+1),...,\calo(N+D+1)$ so that $\ttt_k(W) = 0 \forall k \geq 1$ and $\ttt_0(W)=1$. Clearly, $\ttt_k(\alpha_l(V_i \otimes W)) = \ttt_k(\alpha_l(V_i))$ is killed by $I(k,r,l)$. \\

\end{proof}

Remark: Originally however, the hope was for a stronger result saying that for fixed $l$ and $r$ , there exists a $k$ satisfying $I(k,r,l) \subsetneq I(k,r-1,l)$ . In fact, there was the hope of being able to show that $I(2r, r, l) \subsetneq I(2r, r-1, l)$. This would have shown that there is no morphism $f:G(r, 2r) \rar G(r-1, M)$ so that $f^*(\alpha_l(Q')) = \alpha_l(Q)$. We have so far been unable to do this in general. However, we have found (by means of a computer program) that $I(6,3,2) \subsetneq I(6,2,2)$ thus proving that if $f:G(3,6) \rar G(2,M) $ is a morphism, then $f^*(\alpha_2(Q')) = C {\alpha_1(Q)}^2$. This we do by showing that $\oplus_{|\lambda|=3} \enn(V_{\lambda}) $ contains an irreducible representation $V_{\mu}$ of $S_6$ not contained in $\oplus_{|\lambda| \leq 2} \enn(V_{\lambda})$, and that if $\pi_{\mu}$ denotes the projection from $KS_k$ to $\enn(V_{\mu})$, then $\pi_{\mu} * {\tilde{e}^{(2)}}_6 * \tau_6 \neq 0$. This is achieved using a Mathematica program. \\

\section{ Proof of Theorem 3}

 \subsection{Computing some cup products in $\oplus \Hm^k(G(r,n), {\Omega}^{\otimes k})$}

In this section, we will show how one computes the cup product of two elements $X_k \in \hhh_K({\wedge}^k \Omega, {\Omega}^{\otimes k}) \subset \Hm^k(G(r,n), {\Omega}^{\otimes k})$ and $Y_l \in \hhh_K( {\wedge}^l \Omega, {\Omega}^{\otimes l}) \subset \Hm^l(G(r,n),{\Omega}^{\otimes l})$, where $X_k = (\gamma_k \otimes \id) \circ i_k \in \enn_K({Q^*}^{\otimes k}) \otimes \enn_K(S^{\otimes k})$ and $Y_l = (\delta_l \otimes \id) \circ i_l \in \enn_K({Q^*}^{\otimes l}) \otimes \enn_K(S^{\otimes l})$, where $i_k$ and $i_l$ are the standard inclusions ${\wedge}^k \Omega \rar {\Omega}^{\otimes k}$ and ${\wedge}^l \Omega \rar {\Omega}^{\otimes l}$ respectively, and where $\enn_K(\Omega^{\otimes *})$ is identified with  $\enn_K({Q^*}^{\otimes *}) \otimes \enn_K(S^{\otimes *})$. In this situation we have the following commutative diagrams: (we are following the conventions of Section 4.)\\

\[
 \begin{CD}
0 @> >> \Omega^{\otimes k}@> >> Z_1 @> >> ....@> >>Z_k @> >> K @> >> 0\\
@. @A{\bar{\gamma}_k}AA   @A{\theta_1}AA  @AAA @AAA  @AA{\id}A \\
...@> >>\Omega^{\otimes k} \otimes \sym^* \Omega @> \alpha_k>>\Omega^{\otimes k-1} \otimes \sym^* \Omega@> >> ....@> >> \sym^* \Omega @> >> K @> >>0 \\
\end{CD}
\] \\

\[
 \begin{CD}
0 @> >> \Omega^{\otimes l}@> >> W_1 @> >> ....@> >>W_l @> >> K @> >> 0\\
@. @A{\bar{\delta}_l}AA   @A{\theta_2}AA  @AAA @AAA  @AA{\id}A \\
...@> >>\Omega^{\otimes l} \otimes \sym^* \Omega @> \alpha_l>>\Omega^{\otimes l-1} \otimes \sym^* \Omega@> >> ....@> >> \sym^* \Omega @> >> K @> >>0 \\
\end{CD}
\]

The top rows of the two commutative diagrams are exact sequences representing $X_k$ and $Y_l$ respectively. To compute the cup product $X_k \cup Y_l$ we only need to find vertical arrows making all squares in the following diagram commute:

\[
 \begin{CD}
    0 @> >> \Omega^{\otimes k+l}@> >> Z_1 \otimes \Omega^{\otimes l} @> >> ....@> >>Z_k \otimes \Omega^{\otimes l} @> >> W_1 \\
 @. @AAA @AAA @AAA @AAA  @AAA\\
 ...@> >>\Omega^{\otimes k+l} \otimes \sym^* \Omega @> \alpha_{k+l}>>\Omega^{\otimes k+l-1} \otimes \sym^* \Omega@> >> ....@> >> ...@> >>\Omega^{\otimes l} \otimes \sym^* \Omega
\end{CD}
\]

\[
 \begin{CD}
...@> >> Z_k \otimes \Omega^{\otimes l}@> >> W_1 @> >> ....@> >>W_l @> >> K @> >> 0\\
@. @AAA   @A{\theta_2}AA  @AAA @AAA  @AA{\id}A \\
...@> >> \Omega^{\otimes l} \otimes \sym^* \Omega @> \alpha_l>>\Omega^{\otimes l-1} \otimes \sym^* \Omega@> >> ....@> >> \sym^* \Omega @> >> K @> >>0 \\
\end{CD}
\] \\

Note that the diagrams below commute: \\

    \[
 \begin{CD}
    0 @> >> \Omega^{\otimes k+l}@> >> Z_1 \otimes \Omega^{\otimes l} @> >> ....@> >>Z_k \otimes \Omega^{\otimes l} @> >>\Omega^{\otimes l} \\
 @. @A{\bar{\gamma_k \otimes \delta_l}}AA @A{\theta_1 \otimes \delta_l}AA @AAA @AAA  @AA{\delta_l}A\\
 ...@> >>\Omega^{\otimes k+l} \otimes \sym^* \Omega @> \alpha_{k+l}>>\Omega^{\otimes k+l-1} \otimes \sym^* \Omega @> >> ....@> >> \Omega^{\otimes l} \otimes \sym^* \Omega @> >> {\Omega}^{\otimes l}
\end{CD}
\]\\

\[
 \begin{CD}
    Z_k \otimes \Omega^{\otimes l} @> >>\Omega^{\otimes l}\\
   @A{- \otimes \delta_l}AA             @AA{\delta_l}A\\
    \Omega^{\otimes l} \otimes \sym^* \Omega  @> >> \Omega^{\otimes l}
\end{CD}
\] \\

From these diagrams, we deduce the following lemma: \\
\begin{lemma}
 $[(\gamma_k \otimes \id) \circ i_k] \cup [(\delta_l \otimes \id) \circ i_l] = [((\gamma_k \otimes \delta_l) \otimes \id) \circ i_{k+l}]$ The element $(\gamma_k \otimes \delta_l) \in K(S_k \times S_l) \subset K(S_{k+l})$ where $S_k \times S_l$ is embedded in $S_{k+l}$ in the natural way. \\
\end{lemma}

This enables us to compute cup products of the form $\cup_i {\varphi_i}_* \ttt_{k_i}(\alpha_{l_i}(Q))$. \\

\subsection{A certain decomposition of $KS_k$}

Observe that $KS_k = \oplus W_{\lambda}$ where $W_{\lambda}$ is the $K$-span of elements of $S_k$ in the conjugacy class corresponding to the partition $\lambda$. We shall break each of the spaces $W_{\lambda}$ further into a direct sum of $K$-vector spaces in a specific manner. The significance of the new decomposition shall become clear as we proceed. \\
First, let us decompose the congugacy class $C_{(k)}$ which is the congugacy class of the cycle $\tau_k$. Note that $\tau_k =  \sum_{l \geq 2} {\tilde{e}^{(l)}}_k * \tau_k$ and that  ${\tilde{e}^{(l)}}_k {\tilde{e}^{(l')}}_k = \delta_{ll'}{\tilde{e}^{(l)}}_k$. Define operators $\Pi_l$ on $C_{(k)}$ by $\Pi_l(\beta \tau_k \beta^{-1}) = \beta * ( {\tilde{e}^{(l)}}_k * \tau_k)$ for $\beta \in S_k$ and extend this by linearity to $C_{(k)}$. Note that $\sum_{l \geq 2} \Pi_l(\beta * \tau_k) = \beta * \tau_k$. First, we need to check that we actually have a welldefined operator here. It suffices to show that if $ \beta , \gamma \in S_k$ with $\beta * \tau_k = \gamma * \tau_k$ then $\Pi_l( \beta * \tau_k) = \Pi_l (\gamma * \tau_k)$. In other words, we need to show that $\beta * ( {\tilde{e}^{(l)}}_k * \tau_k) = \gamma * ( {\tilde{e}^{(l)}}_k * \tau_k)$ which is eqvivalent to showing that $(\beta^{-1} \gamma) * ( {\tilde{e}^{(l)}}_k * \tau_k) =  {\tilde{e}^{(l)}}_k * \tau_k$. But $\beta * \tau_k = \gamma * \tau_k$ iff $\beta^{-1} \gamma = {\tau_k}^s $ for some $s$. Therefore, that $\Pi_l$ is well defined follows from the following lemma: \\

\begin{lemma}
$ {\tau_k}^s *(  {\tilde{e}^{(l)}}_k * \tau_k) =  {\tilde{e}^{(l)}}_k * \tau_k$ for any integer $s$.
\end{lemma}

\begin{proof}
This really follows from the fact that for any smooth scheme $X$, and for any vector bundle $V$ on $X$, $\sn(\tau_k) {\tau_k}_* \ttt_k(V) = \ttt_k(V)$. After all, $\sn(\tau_k) {\tau_k}_* {\theta_V}^k = {\theta_V}^k $ (by the properties of the cup product). Hence, $tr_* \varphi_* \sn(\tau_k) {\tau_k}_* {\theta_V}^k = tr_* \varphi_* {\theta_V}^k$ . The right hand side of this equation is $\ttt_k(V)$ by definition. The left hand side is   $\sn(\tau_k) {\tau_k}_* \ttt_k(V)$ since $tr \circ \varphi \circ \tau_k = \tau_k \circ tr \circ \varphi$. \\
This tells us that $\sn({\tau_k}^s) {{\tau_k}^s}_* \ttt_k(V) = \ttt_k(V)$.  To finishthe proof of the lemma , we observe that  ${\tau_k}^s *(  {\tilde{e}^{(l)}}_k * \tau_k) = \sn({\tau_k}^s) {{\tau_k}^s}_* \ttt_k(\alpha_l(Q))$. \\
\end{proof}

The other detail to be verified is the fact that the operators $\Pi_l$ are mutually orthogonal projections. For this, we see that $\Pi_l(\beta *\tau_k) = \beta * (  {\tilde{e}^{(l)}}_k * \tau_k) = (\beta {\tilde{e}^{(l)}}_k) * \tau_k \implies \Pi_l \circ \Pi_m (\beta * \tau_k) = (\beta {\tilde{e}^{(m)}}_k {\tilde{e}^{(l)}}_k)* \tau_k  = (\beta \delta_{lm} {\tilde{e}^{(l)}}_k) *\tau_k  $. We therefore, have a direct sum decomposition  $W_{(k)} = \oplus_{l \geq 2} \Pi_l(W_{(k)})$. \\
We now proceed to breakup $W_{\lambda}$ into a direct sum of $K$-vector spaces in an analogous manner. Note that $C_{\lambda}$ is the conjugacy class of $\tau_{\lambda} := \tau_{\lambda_1}  \tau_{\lambda_2} .... \tau_{\lambda_s}$ where the partition $\lambda$ is givan by $\lambda: k = \lambda_1 + .. + \lambda_s$ , the $\lambda_i$ 's arranged in decreasing order, and where $\tau_{\lambda_i}$ is the cycle $(\lambda_1 + ... + \lambda_{i}, \lambda_1 + ... + \lambda_{i}-1,...,\lambda_1 + ... + \lambda_{i-1})$ which is after all the cycle $\tau_{\lambda_i}$ embedded in $S_k$ under the composition $S_{\lambda_i} \subset S_{\lambda_1} \times ... \times S_{\lambda_s} \subset S_k$. Call the map $S_{\lambda_1} \times ... \times S_{\lambda_s} \subset S_k$ as $\varphi$. Note that $\varphi$ extends to a $K$-algebra homomorphism $\varphi: K(S_{\lambda_1} \times ... \times S_{\lambda_s}) \rar K(S_k)$. Identify $  K(S_{\lambda_1}) \otimes ... \otimes K(S_{\lambda_s})$ with $ K(S_{\lambda_1} \times ... \times S_{\lambda_s})$ and consider $ ({\tilde{e}^{(l_1)}}_{\lambda_1} \otimes .... \otimes {\tilde{e}^{(l_s)}}_{\lambda_s})*\tau_{\lambda}$. By this we are looking at ${\tilde{e}^{(l_1)}}_{\lambda_1} \otimes .... \otimes {\tilde{e}^{(l_s)}}_{\lambda_s}$ as an element of $KS_k$ through the homomorphism $\varphi$. We now make the following observations: \\

Observation 1: The elements ${\tilde{e}^{(l_1)}}_{\lambda_1} \otimes .... \otimes {\tilde{e}^{(l_s)}}_{\lambda_s}$ are mutually orthogonal idempotents in $K(S_k)$ adding up to $\id$. This follows from the fact that the above statement is true in $K(S_{\lambda_1} \times ... \times S_{\lambda_s})$. \\

Observation 2: As $\tau_{\lambda} = \tau_{\lambda_1} \otimes ... \otimes \tau_{\lambda_s}$, $({\tilde{e}^{(l_1)}}_{\lambda_1} \otimes .... \otimes {\tilde{e}^{(l_s)}}_{\lambda_s})*\tau_{\lambda} = ( {\tilde{e}^{(l_1)}}_{\lambda_1} * \tau_{\lambda_1}) \otimes ... \otimes ({\tilde{e}^{(l_s)}}_{\lambda_s} * \tau_{\lambda_s})$. It follows that if for some $i$ , $\lambda_i \geq 2$ and $l_i =1$, then $({\tilde{e}^{(l_1)}}_{\lambda_1} \otimes .... \otimes {\tilde{e}^{(l_s)}}_{\lambda_s})*\tau_{\lambda} =0$. \\

Observation 3: Let $ {\tilde{e}^{(l)}}_{\lambda} := \sum_{l_1 +...+l_s =l} {\tilde{e}^{(l_1)}}_{\lambda_1} \otimes .... \otimes {\tilde{e}^{(l_s)}}_{\lambda_s}$. Then ${\tilde{e}^{(l)}}_{\lambda}$ is an idempotent with ${\tilde{e}^{(l)}}_{\lambda} . ( {\tilde{e}^{(l_1)}}_{\lambda_1} \otimes .... \otimes {\tilde{e}^{(l_s)}}_{\lambda_s}) = ({\tilde{e}^{(l_1)}}_{\lambda_1} \otimes .... \otimes {\tilde{e}^{(l_s)}}_{\lambda_s})$ if $l_1+...+l_s =l$ and ${\tilde{e}^{(l)}}_{\lambda} . ( {\tilde{e}^{(l_1)}}_{\lambda_1} \otimes .... \otimes {\tilde{e}^{(l_s)}}_{\lambda_s}) = 0$ otherwise. \\

Let $\Pi_l$ be defined by $\Pi_l( \beta * \tau_{\lambda}) = (\beta {\tilde{e}^{(l)}}_{\lambda}) * \tau_{\lambda}$. We then have :\\

\begin{lemma}
The $\Pi_l$ are well-defined mutually orthogonal projection operators on $W_{\lambda}$.
\end{lemma}

\begin{proof}
Note that it suffices to show that if $\gamma$ is a permutation in the stabiliser of $\tau_{\lambda}$ under conjugation, then $\gamma*({\tilde{e}^{(l)}}_{\lambda} * \tau_{\lambda}) = {\tilde{e}^{(l)}}_{\lambda} * \tau_{\lambda}$. Note that if $\gamma$ stabilises $\tau_{\lambda}$ under conjugation, then $\gamma$ is of the form $\zeta ( {\tau_{\lambda_1}}^{r_1} \otimes .... \otimes {\tau_{\lambda_s}}^{r_s})$ where $\zeta$ permutes blocks of equal lengths among $[1,...,\lambda_1],[\lambda_1+1,...,\lambda_1+ \lambda_2],....,[\lambda_1+...+\lambda_{s-1}+1,...,k]$ while preserving order within such blocks. Now we need to show that $\gamma*({\tilde{e}^{(l)}}_{\lambda} * \tau_{\lambda}) = {\tilde{e}^{(l)}}_{\lambda} * \tau_{\lambda}$. Observe that $( {\tau_{\lambda_1}}^{r_1} \otimes .... \otimes {\tau_{\lambda_s}}^{r_s}) * ({\tilde{e}^{(l_1)}}_{\lambda_1} \otimes .... \otimes {\tilde{e}^{(l_s)}}_{\lambda_s}) * \tau_{\lambda} = ({\tau_{\lambda_1}}^{r_1} * {\tilde{e}^{(l_1)}}_{\lambda_1} * \tau_{\lambda_1} ) \otimes .... \otimes ({\tau_{\lambda_s}}^{r_s} * {\tilde{e}^{(l_s)}}_{\lambda_s} * \tau_{\lambda_s} ) = ({\tilde{e}^{(l_1)}}_{\lambda_1} \otimes .... \otimes {\tilde{e}^{(l_s)}}_{\lambda_s}) * \tau_{\lambda}$  (the last equality by Lemma 19). So, we only need to show that $\zeta * {\tilde{e}^{(l)}}_{\lambda} * \tau_{\lambda} ={\tilde{e}^{(l)}}_{\lambda} * \tau_{\lambda}$. But this is true since $\zeta$ induces a permutation $\zeta'$ of $1,2,..,s$ and we see that $\zeta .({\tilde{e}^{(l_1)}}_{\lambda_1} \otimes .... \otimes {\tilde{e}^{(l_s)}}_{\lambda_s}) = ({{\tilde{e}}^{(l_{\zeta'(1)})}}_{\lambda_{\zeta'(1)}} \otimes .... \otimes {{\tilde{e}}^{(l_{\zeta'(s)})}}_{\lambda_{\zeta'(s)}})$  \\
\end{proof}

Observation 4: It now follows from this and the fact that the $\Pi_l$ are mutually orthogonal idempotents adding upto $\id$ that $W_{\lambda} = \oplus \Pi_l(W_{\lambda})$. Also, Observation 2. tells us that $\Pi_1(W_{\lambda}) = 0$ and that  $\Pi_2(W_{\lambda}) = 0$ if $\lambda \neq (k)$. Therefore, this direct sum decomposition runs over $l \geq 2$. Combining this with the decomposition $KS_k = \oplus_{\lambda} W_{\lambda}$, we see that $KS_k = \oplus_{\lambda} \oplus_{l \geq 2} \Pi_l(W_{\lambda}) = \oplus_{l \geq 2} \Pi_l(KS_k) $.\\

 Remark: We now have to see what this means in terms of our problem. Lemma 18 tells us that $\ttt_{\lambda_1}(\alpha_{l_1}(Q)) \cup ... \cup \ttt_{\lambda_s}(\alpha_{l_s}(Q)) = [({\tilde{e}^{(l)}}_{\lambda} * \tau_{\lambda}) \otimes \id] \circ i$ as an element of $\hhh_K({\wedge}^k \Omega, \Omega^{\otimes k})$. In the notation of Sections 4 and 5, $\sn(\beta) {\beta^{-1}}_* \ttt_{\lambda_1}(\alpha_{l_1}(Q)) \cup ... \cup \ttt_{\lambda_s}(\alpha_{l_s}(Q)) = (\Pi_l(\beta * \tau_k) \otimes \id) \circ i$ in $\hhh_K({\wedge}^k \Omega, \Omega^{\otimes k})$. Our observations tell us the following: After all $\ttt_{\lambda_1}(\alpha_{l_1}({\psi}^pQ)) \cup ... \cup \ttt_{\lambda_s}(\alpha_{l_s}({\psi}^pQ)) = p^l \ttt_{\lambda_1}(\alpha_{l_1}(Q)) \cup ... \cup \ttt_{\lambda_s}(\alpha_{l_s}(Q))$ where $\sum l_i =l$.Thus the space spanned by $\beta_* \ttt_{\lambda_1}(\alpha_{l_1}(Q)) \cup ... \cup \ttt_{\lambda_s}(\alpha_{l_s}(Q))$, which is $\Pi_l(W_{\lambda})$ can be thought of as a $'\CH^l$ (Adams weight $l$) space in $W_{\lambda}$. We have, in this section shown that $KS_k$ decomposes into a direct sum of these ``Adams weight $l$'' spaces that are stable under conjugation by elements of $S_k$ by their very construction. \\
This states that we cannot hope to find some expresstion $X(Q)=\beta_* \ttt_{\lambda_1}(\alpha_{l_1}(Q)) \cup ... \cup \ttt_{\lambda_s}(\alpha_{l_s}(Q))$ of Adams weight $l:= \sum l_i$ equal to some other expression $Y(Q)$ of a different Adams weight in $KS_k$. If that were the case, then $(X-Y)(Q) = 0$ and $(X-Y)({\psi}^pQ) \neq 0$ with the expressions $X()$ and $Y()$ being functorial with respect to pullbacks. That would have solved our problem straightaway. \\
However, even though there is no linear dependence relation among expressions of different Adams weights is $KS_k$, we will show that such a relation exists in $\hhh_K({\wedge}^k \Omega, \Omega^{\otimes k})$ if we are looking at a Grassmannian $G(r,n)$. That will be done by showing that the decomposition of $KS_k$ with respect to Adams spaces does not ``behave well'' with respect to the projection to the space $\oplus_{|\lambda| \leq r} \enn(V_{\lambda})$.\\

\subsection{A linear dependence relation between functors of different Adams weights}

First, we observe that if $V$ is a vector space with $V= V_1 \oplus V_2$ and also $V = \oplus W_i$, with $p_i$ being the projections to $V_i$ and $\pi_i$ being the projections to $W_i$, then $\dim \text{ } p_1(W_1)) +...+ \dim \text{ } p_1(W_m) \geq \dim \text{ }V_1$. Suppose that equality holds. Then $\dim \text{ } p_1(W_i) = \dim \text{ } W_i - \dim \text{ } W_i \cap V_2  \implies \dim \text{ } W_1 \cap V_2 +... +\dim \text{ } W_m \cap V_2 =\dim \text{ } V_2 $. From this, we see that $\pi_i(V_2) = W_i \cap V_2$ for all $i \in \{1,2,..,,\}$. In particular, if $\pi_i(V_2) \neq W_i \cap V_2$, then $\dim \text{ } p_1(W_1)) +...+ \dim \text{ } p_1(W_m) > \dim \text{ }V_1$. \\
Having said this, we will prove that for $V = KS_{2r}$ ( $V = V_1 \oplus V_2$ where $V_1 = \oplus_{|\lambda| \leq r} \enn(V_{\lambda})$ and $V_2 = \oplus_{|\lambda| > r} \enn(V_{\lambda})$ also $V = \oplus_{l \geq 2} \Pi_l(V)$) $\Pi_2(V_2) \neq \Pi_2(V) \cap V_2$ for some $r$ at least, we we shall later specify. This will prove that $\sum_{l \geq 2}  \dim \text{ } P_r(\Pi_l(V)) > \dim \text{ }V_1$ for these $r$. (Note that $P_r$ here is $p_1$ for this situation). The observations in the previous subsection tell us that $\Pi_2(V) = \Pi_2(W_{(2r)}$. Any element in this space is a linear combination of conjugates of $\tau_2r$. It follows that if such a linear combination is nonzero in $\enn(V_{\lambda})$ it is also nonzero as an element of $\enn(V_{\bar{\lambda}})$, where ${\bar{\lambda}}$ is the partition conjugate to $\lambda$. Thus $\Pi_{(2r)}(V) \cap V_2 =0$. It therefore , suffices to prove that $\Pi_2(V_2) \neq 0$. \\

\begin{lemma}
 It suffices to show that $\Pi_2 ( (1 \text{ } 2r) \sum_{g \in C_{\mu}} g) \neq 0$ where $(1 \text{ }  2r)$ is the transposition interchanging $1$ with $2r$ and $\mu$ is a partition among $\{(2r-1,1),...,(r,r)\}$.
\end{lemma}

\begin{proof} Consider the matrix $M = (\chi_{\lambda}(C_{\mu}))$ where $\lambda$ runs over all partitions of $2r$ that satisfy $\lambda \geq (r,r)$ (There is a lexicographic ordering among the partitions,enabling one to compare them), and $\mu \in \{(2r-1,1),...,(r,r)\}$. Note that is $\lambda$ is such a partition and $\lambda \neq (r,r)$ then $\lambda_1 \geq r+1$. We claim that $M$ is of rank $r$. To prove this, it suffices to show that $N$ is of rank $r$ where $N = (\psi_{\lambda}(C_{\mu}))$, where $\psi_{\lambda} =  {Ind_{S_{\lambda}}}^{S_{2r}} (triv) = \chi_{\lambda} + \sum_{\mu > \lambda} K_{\mu \lambda} \chi_{\mu}$. However,   $\psi_{\lambda}(C_{\mu}) = \frac{1}{|C_{\mu}|}[S_{2r}:S_{\lambda}]|C_{\mu} \cap S_{\lambda}|$. Therefore, $ \psi_{\lambda}(C_{\mu}) =0$ if $\mu > \lambda$. This lexicographic order is a total order. Consider the restriction of $N$ to the rows given by the partitions in $\{(2r-1,1),...,(r,r)\}$. This restriction of $N$ is then a lower triangular matrix with nonzero diagonal entries if the rows are arranged in the correct order (since $ \psi_{\lambda}(C_{\lambda})  \neq 0$).It follows that $N$, and hence $M$ are matrices of rank $r$. We further claim that if we restrict $M$ to rows corresponding to $\lambda > (r,r)$, we still get a matrix of rank $r$. To see this, we need to show that for some scalars $a_{\lambda}$ , $\chi_{(r,r)}(C_{\mu}) = \sum_{\lambda > (r,r)} a_{\lambda}\chi_{\lambda}(C_{\mu})$ for all $\mu \in \{(2r-1,1),...,(r,r)\}$. For this, it is enough to show that $\psi_{(r,r)}(C_{\mu}) = \sum_{\lambda > (r,r)} b_{\lambda}\psi_{\lambda}(C_{\mu})$ for all $\mu \in \{(2r-1,1),...,(r,r)\}$, for some scalars $b_{\lambda}$. In fact, we claim that there are scalars $b_i , 0 \leq i \leq r-1$, so that $\psi_{(r,r)}(C_{\mu}) = \sum_{0 \leq i \leq r-1} b_i \psi_{(2r-i,i)}(C_{\mu})$. Note that $|C_{(2r-s,s)} \cap S_{(2r-t,t)}| = 0$ if $s \neq t$ and both are nonzero. Also note that $\psi_{(2r)}(C_{(r,r)}) \neq 0$. Thus the vector $\psi_{(2r)}(C_{\mu}) , \mu \in \{(2r-1,1),...,(r,r)\}$ is given by $(a_1,..,a_r)$, where$a_r \neq 0$. The vector $\psi_{(2r-s,s)}(C_{\mu}) , \mu \in \{(2r-1,1),...,(r,r)\}$  is given by $(0,..,o,d_s,...,0)$ , $d_s \neq 0$ for $1 \leq s \leq r-1$. This $\psi_{(2r)}(C_{\mu}) - \sum \frac{a_s}{d_s} \psi_{(2r-s,s)}(C_{\mu}) = (0,..,0,a_r)$ which is a nonzero multiple of $\psi_{(r,r)}(C_{\mu})$. This shows that the matrix $M = \chi_{\lambda}(C_{\mu})$ where $\lambda > (r,r)$ and $\mu \in  \{(2r-1,1),...,(r,r)\}$ is of rank $r$. Since $\chi_{\bar{\lambda}} = \chi_{\lambda}.\sn$, and $|\bar{\lambda}| \geq r+1$ iff $\lambda > (r,r)$, the matrix $M' = \chi_{\lambda}(C_{\mu})$ where $|\bar{\lambda}| \geq r+1 $ and $\mu \in  \{(2r-1,1),...,(r,r)\}$ is obtained from $M$ by multiplying some columns by $-1$ and is thus of rank $r$. \\
Now suppose that $\Pi_2( (1 \text{ } 2r) \sum_{g \in C_{(2r-s,s)}} g) \neq 0$ for somw $1 \leq s \leq r$ . Since $M'$ is of rank $r$, we can find a linear combination of rows of $M'$ that gives us the vector $e_s$ i.e, $\sum_{|\lambda| > r+1} a_{\lambda} \chi_{\lambda}(C_{\mu}) = 0$ if $\mu \neq (2r-s,s)$ and  $\sum_{|\lambda| > r+1} a_{\lambda} \chi_{\lambda}(C_{\mu}) = 1$ if $\mu = (2r-s,s)$. So, $\Pi_2( (1 \text{ } 2r) (\sum_{g \in S_{2r} ; |\lambda| > r+1} a_{\lambda} \chi_{\lambda}(g)g)) = \Pi_2( (1 \text{ } 2r) \sum_{g \in C_{(2r-s,s)}} g) \neq 0$. The first equality is because only the $2r$ cycles contribute to $\Pi_2(V)$. Note that if $ \sum \chi_{\lambda}(g)g \in \enn(V_{\lambda})$ it follows that $(\sum_{g \in S_{2r} ; |\lambda| > r+1} a_{\lambda} \chi_{\lambda}(g)g) \in V_2$ and thus \\  $(1 \text{ } 2r) (\sum_{g \in S_{2r} ; |\lambda| > r+1} a_{\lambda} \chi_{\lambda}(g)g) \in V_2$ thus proving that $\Pi_2(V_2) \neq 0$. \\
\end{proof}

\begin{lemma}
For some $s$, $1 \leq s \leq r$, we have  $\Pi_2 ( (1 \text{ } 2r) \sum_{g \in C_{(2r-s,s)}} g) \neq 0$
\end{lemma}

\begin{proof}
Every $2r$ cycle that arises in $(1 \text{ } 2r) \sum_{g \in C_{(2r-s,s)}} g$ arises with coefficient $1$. We therefore need to identify the $2r$ cycles that do arise. They are those of the form $(1 \text{ } a_2 \text{ } .. \text{ } a_s \text{ } 2r \text{ } a_{s+2} \text{ } ....)$  or $(1 \text{ } a_2 \text{ }  .... \text{ } a_{2r-s} \text{ } 2r ....)$. We note that
 $$ (1 \text{ } 2r) \sum_{g \in C_{(2r-s,s)}} g$$
 $$ = \sum_{\alpha \in S_{2r-1} \text{ fixing } 1 \text{ and } 2r} \alpha * (2r \text{ } 2r-s \text{ } 2r-s-1 .. .. 1 \text{ }2r-1\text{ } 2r-2 ...\text{ } 2r-s+1)$$
 $$ + \alpha* (2r \text{ }s \text{ } s-1 ... \text{ } 1 \text{ } 2r-1 ..\text{ } s+1)$$

\begin{eqnarray*}
& = & \sum_{\alpha \in S_{2r-1} \text{ fixing } 1 \text{ and } 2r} \alpha * [ {\tau^{s-1}}_{2r-1} + {\tau^{2r-s-1}}_{2r-1}]* \tau_{2r}\\
& = & [{\tau^{-1}}_{2r-1} [\sum_{\beta \in S_{2r-1} \text{ fixing } 2r-1 \text{ and } 2r} \beta] \tau_{2r-1}] * [ {\tau^{s-1}}_{2r-1} + {\tau^{2r-s-1}}_{2r-1}]* \tau_{2r}\\
& = & [{\tau^{-1}}_{2r-1}  \sum_{\beta \in S_{2r-1} \text{ fixing } 2r-1 \text{ and } 2r} \beta ] * [{\tau^{s}}_{2r-1} + {\tau^{2r-s}}_{2r-1}]* \tau_{2r}
\end{eqnarray*}
Therefore, $$ \Pi_2( (1 \text{ } 2r) \sum_{g \in C_{(2r-s,s)}} g) = [{\tau^{-1}}_{2r-1}  \sum_{\beta \in S_{2r-1} \text{ fixing } 2r-1 \text{ and } 2r} \beta ] * [{\tau^{s}}_{2r-1} + {\tau^{2r-s}}_{2r-1}]* [{\tilde{e}^{(2)}}_{2r} * \tau_{2r}]$$ $$ =  [{\tau^{-1}}_{2r-1}  \sum_{\beta \in S_{2r-1} \text{ fixing } 2r-1 \text{ and } 2r} \beta ] * [{\tau^{s}}^{2r-1} + {\tau^{2r-s}}_{2r-1}]* [{\tilde{e}^{(1)}}_{2r-1} * \tau_{2r}]$$ , the last equality by Lemma 16. \\
It therefore, suffices to show that $  [{\tau^{-1}}_{2r-1}  \sum \beta ] [{\tau^{s}}^{2r-1} + {\tau^{2r-s}}_{2r-1}] [{\tilde{e}^{(1)}}_{2r-1}] \neq 0$ for some $s$, $1 \leq s \leq r$. It therefore, suffices to show that $ W_s :=[\sum \beta ] [{\tau^{s}}^{2r-1} + {\tau^{2r-s}}_{2r-1}] [{\tilde{e}^{(1)}}_{2r-1}] \neq 0$ for some $s$. Consider a vector space $V$ of finite dimension, and let $u$ and $v$ be two basis vectors of $V$. We will show that the right action of $W_s$ on $u^{\otimes 2r-2} \otimes v$ is nonzero. Note that $\frac{1}{(2r-2)!} (u^{\otimes 2r-2} \otimes v) W_s = (u^{\otimes 2r-2} \otimes v) ({\tau^{s}}^{2r-1} + {\tau^{2r-s}}_{2r-1}) {\tilde{e}^{(1)}}_{2r-1} = (u^{\otimes s-1} \otimes v \otimes u^{\otimes 2r-1-s}  + u^{\otimes 2r-1-s} \otimes v \otimes u^{\otimes s-1}) {\tilde{e}^{(1)}}_{2r-1} $.Therefore it is enough to show that  $(u^{\otimes s-1} \otimes v \otimes u^{\otimes 2r-1-s}  + u^{\otimes 2r-1-s} \otimes v \otimes u^{\otimes s-1}) {\tilde{e}^{(1)}}_{2r-1} \neq 0$ for some $s$. For this, we note that $ad(u)^{2r-2}(v) = (l_u - r_u)^{2r-2} (v) = \sum_i  \binom{2r-2}{i} u^{\otimes i} \otimes v \otimes u^{2r-2-i}$ The idempotent $  {\tilde{e}^{(1)}}_{2r-1}$ acts as the identity on this vector, which is a linear combination of  $(u^{\otimes s-1} \otimes v \otimes u^{\otimes 2r-1-s}  + u^{\otimes 2r-1-s} \otimes v \otimes u^{\otimes s-1})$ where $s$ runs from $1$ to $r$. \\
\end{proof}

We have already noted that $\Pi_l( KS_{2r})$ is the space of the elements o $KS_{2r}$ representing  those $\gamma_* \ttt_{\lambda_1}(\alpha_{l_1}(Q)) \cup ... \cup \ttt_{\lambda_s}(\alpha_{l_s}(Q))$ whose Adams weight is $l$ , i.e, $l_1+..+l_s =l$. If $V_1 = \sum_{| \lambda | \leq r} \enn(V_{\lambda})$ and $p_1$ denotes the projection from $KS_{2r}$ to $V_1$, then any element in $p_1( \Pi_l(KS_{2r}))$ is an element of $\hhh_K( {\wedge}^k \Omega, \Omega^{\otimes k})$ that can be expressed as a linear combination of expressions of the form $\gamma_* \ttt_{\lambda_1}(\alpha_{l_1}(Q)) \cup ... \cup \ttt_{\lambda_s}(\alpha_{l_s}(Q))$ of Adams weight $l$. That $\sum_l \dim(p_1(\Pi_l(KS_{2r}))) > \dim(V_1)$ tells us that there is a linear dependence relation among nonzero vectors of different Adams weight, $\sum_l v_l(Q) =0$, where each $v_l$ is a linear combination of expressions of the form $\gamma_* \ttt_{\lambda_1}(\alpha_{l_1}(. )) \cup ... \cup \ttt_{\lambda_s}(\alpha_{l_s}(.))$ such that $l_1 +...+l_s =l$. Chose such a linear dependence relation of minimal length, and observe that $v_l( \psi^p Q) = p^l v_l(Q)$.Also observe that the $v_l$ are functorial under pullback. This leads us to the following result :\\

\begin{theorem}

Let $Q$ denote the universal quotient bundle of a Grassmannian $G(r,n)$ , $r \geq 2$ , $n \geq 2r+1$. Then, for all $p \geq 2$, $[ \psi^p Q]$ is not equal, in K-theory to the class of a genuine vector bundle.

\end{theorem}

\begin{proof}
Suppose that $ [ \psi^p Q] =[Y]$ for some genuine vector bundle $Y$. Then $Y$ is of rank $r$, and for all sufficiently large $m$ , $Y \otimes \calo(m)$ is a quotient of ${\calo_G}^s$ for some $s$. It follows that $Y \otimes \calo(m) = f^*Q$ for some morphism $f: G(r,n) \rar G(r,n')$, where $Q$ is the universal quotient bundle of $G(r,n')$. Note that we have $\sum_l v_l(Q) = 0 \implies \sum_l v_l(Y \otimes \calo(m) = 0$ for all sufficiently large $m$. Note that $\oplus \ttt_k(\calo(m)) = exp( \ttt_1(\alpha_1(\calo(1))))$. Therefore, $\ttt_{\lambda_i}(\alpha_{l_i}(Y \otimes \calo(m))) = \ttt_{\lambda_s}(\alpha_{l_s}(Y) + m \alpha_{l_s-1}(Y) \alpha_1(\calo(1)) + ....)$. Therefore, $v_k(Y \otimes \calo(m)) = v_k(Y) + m.v_{k1}(Y) + ... +m^s v_{ks}(Y)$ for all $k$. It follows from the invertibility of a Vandermonde determinant that  $\sum_l v_l(Y \otimes \calo(m)) = 0$ for all sufficiently large $m$ gives us  $\sum_l v_l(Y) =0$. Thus $\sum_l p^l v_l(Q) =0$ as well. As $p \geq 2$, this contradicts the fact that the linear dependence relation $\sum_l v_l(Q) =0$ was chosen to be of minimal length.\\
\end{proof}

\section{A formula for $\ttt_k(V)$ in terms of $\ch_l(V)$}

In this section, we prove that if $X$ is a projective variety, then $\ttt_k(V)$ can be written in terms of $\ch_l(V) , l \leq k$ and something intrinsic to the given variety. To be more precise, $\ttt_k(V) = \sum_{l \leq k} x_l(V)$, where $x_l(V)$ is obtained from $\ch_l(V)$ by applying finitely many of the following operations: \\
1. Composition in the Yoneda sense with $\psi_{X/S} \in  \eee^1( \Omega_{X/S} , \sym^2 \Omega_{X/S})$ where $ \psi_{X/S} = p_{1*}(0 \rightarrow {\mathcal I}^2/{\mathcal I}^3 \rightarrow {\mathcal I}/{\mathcal I}^3 \rightarrow  {\mathcal I}/ {\mathcal I}^2 \rightarrow 0) =: 0 \rightarrow \sym^2{\Omega} \rightarrow { E} \rightarrow {\Omega} \rightarrow 0$. (Here, $\mathcal I$ denotes the sheaf of ideals defining the diagonal in $X \times_S X$ and $p_1$ is the first pojection from $X \times_S X$ to $X$).\\
2. Actions of elements of $KS_k$ on $\eee^k(\calo_X , \Omega^{\otimes k})$. \\

\subsection{Preliminary results}

 \begin{lemma}
$\tilde{\ttt_3}(V) = \tilde{\ch_3}(V) - p_3 \circ (\id_V \otimes \psi_{X/S} \otimes \id_{\Omega}) \circ \tilde{\ch_2}(V) + Y$ where $Y$ vanishes on taking the trace. Here, $p_3 : \Omega^{\otimes 3} \rar \sym^3 \Omega$ is the standard projection.
\end{lemma}

\begin{proof}

Let $\psi_V$ denote the element in $\eee^1(V \otimes \Omega_{X/S} ,V \otimes \sym^2 \Omega_{X/S})$ given by the exact sequence $  p_{1*}(0 \rightarrow {\mathcal I}^2/{\mathcal I}^3 \otimes_{\calo_Y} {p_2}^*V  \rightarrow {\mathcal I}/{\mathcal I}^3 \otimes_{\calo_Y} {p_2}^*V \rightarrow  {\mathcal I}/ {\mathcal I}^2 \otimes_{\calo_Y} {p_2}^*V \rightarrow 0)$. We first prove the following claim:\\
\begin{claim}
$\psi_V \circ \theta_V =0$
\end{claim}

\begin{proof}
For this, we need the fact that $\theta_V$ is given by the exact sequence $  p_{1*}(0 \rightarrow {\mathcal I}/{\mathcal I}^2 \otimes_{\calo_Y} {p_2}^*V  \rightarrow {\calo_Y}/{\mathcal I}^2 \otimes_{\calo_Y} {p_2}^*V \rightarrow  {\calo_Y}/ {\mathcal I} \otimes_{\calo_Y} {p_2}^*V \rightarrow 0)$. Therefore, $\psi_V \circ \theta_V$ is given by the exact sequence $  p_{1*}(0 \rightarrow {\mathcal I}^2/{\mathcal I}^3 \otimes_{\calo_Y} {p_2}^*V  \rightarrow  {\mathcal I}/{\mathcal I}^3 \otimes_{\calo_Y} {p_2}^*V  \rightarrow {\calo_Y}/{\mathcal I}^2 \otimes_{\calo_Y} {p_2}^*V \rightarrow  {\calo_Y}/ {\mathcal I} \otimes_{\calo_Y} {p_2}^*V \rightarrow 0)$. Note that all squares in the diagram below commute: (where $Z = p_2^*V$)\\

\[
\begin{CD}
  0 @> >> {\mathcal I}^2/{\mathcal I}^3 \otimes_{\calo_Y} Z  @> >>  {\mathcal I}/{\mathcal I}^3 \otimes_{\calo_Y} Z  @> >> {\calo_Y}/{\mathcal I}^2 \otimes_{\calo_Y} Z @> >>  {\calo_Y}/ {\mathcal I} \otimes_{\calo_Y} Z @> >>  0\\
@. @AAA @A{\id}AA @AAA @AA{\id}A @. \\
 0 @> >> 0  @> >>  {\mathcal I}/{\mathcal I}^3 \otimes_{\calo_Y} Z  @> >> {\calo_Y}/{\mathcal I}^3 \otimes_{\calo_Y} Z @> >>  {\calo_Y}/ {\mathcal I} \otimes_{\calo_Y} Z @> >>  0\\
 \end{CD}
\]\\

This shows that $\psi_V \circ \theta_V$ is induced by an element of $\eee^2(V,0) =0$, and is therefore $0$.
\end{proof}

\begin{claim}
$\psi_V = \id_V \otimes_X \psi_{X/S} + (\id_V \otimes p_2) \circ (\theta_V \otimes_X \id_{\Omega_{X/S}})$
\end{claim}

\begin{proof}
Let $W = {p_2}^*V$ and let $ \otimes$ denote $\otimes_Y$. Then, $\psi_V$ is given by $0 \rar  {\mathcal I}^2 W/{\mathcal I}^3 W  \rar  {\mathcal I} W/{\mathcal I}^3 W \rar  {\mathcal I} W/{\mathcal I}^2 W \rar 0$. $\id_V \otimes_X \psi_{X/S}$ is given by $0 \rar W/{\mathcal I}W \otimes_X  {\mathcal I}^2/{\mathcal I}^3   \rar W/{\mathcal I}W \otimes_X  {\mathcal I}/{\mathcal I}^3  \rar W/{\mathcal I}W \otimes_X  {\mathcal I}/{\mathcal I}^2  \rar 0$. $(\id_V \otimes p_2) \circ (\theta_V \otimes_X \id_{\Omega_{X/S}})$, is given by an exact sequence  $0 \rar W/{\mathcal I}W \otimes_X  {\mathcal I}^2/{\mathcal I}^3   \rar E  \rar W/{\mathcal I}W \otimes_X  {\mathcal I}/{\mathcal I}^2  \rar 0$ . The following diagram commuts:\\

\[
\begin{CD}
  0 @> >> {\mathcal I}W/{\mathcal I}^2 W \otimes_X  {\mathcal I}/{\mathcal I}^2   @> >>  W/{\mathcal I}^2 W \otimes_X  {\mathcal I}/{\mathcal I}^2  @> >> W/{\mathcal I}W \otimes_X  {\mathcal I}/{\mathcal I}^2  @> >>  0\\
@. @Vp_2VV @VVV @VV{\id}V @.\\
0 @> >> W/{\mathcal I} W \otimes_X  {\mathcal I}^2/{\mathcal I}^3   @> >> E  @> >> W/{\mathcal I}W \otimes_X  {\mathcal I}/{\mathcal I}^2  @> >>  0\\
  \end{CD}
\]

To prove the desired result consider the following diagram all of whose squares commute , and whose columns are also short exact sequences\\

\[
\begin{CD}
   0 @> >> {\mathcal I}W/{\mathcal I}^2 W \otimes_X  {\mathcal I}^2/{\mathcal I}^3   @> >>  W/{\mathcal I}^2 W \otimes_X  {\mathcal I}^2/{\mathcal I}^3  @> >> W/{\mathcal I}W \otimes_X  {\mathcal I}^2/{\mathcal I}^3  @> >>  0\\
 @. @VVV @VVV @VVV @. \\
 0 @> >> {\mathcal I}W/{\mathcal I}^2 W \otimes_X  {\mathcal I}/{\mathcal I}^3   @> >>  W/{\mathcal I}^2 W \otimes_X  {\mathcal I}/{\mathcal I}^3  @> >> W/{\mathcal I}W \otimes_X  {\mathcal I}/{\mathcal I}^3  @> >>  0\\
@. @VVV @VVV @VVV @. \\
 0 @> >> {\mathcal I}W/{\mathcal I}^2 W \otimes_X  {\mathcal I}/{\mathcal I}^2   @> >>  W/{\mathcal I}^2 W \otimes_X  {\mathcal I}/{\mathcal I}^2  @> >> W/{\mathcal I}W \otimes_X  {\mathcal I}/{\mathcal I}^2  @> >>  0\\
\end{CD}
\]

We then get an exact sequence $0 \rar W/{\mathcal I} W \otimes_X  {\mathcal I}^2/{\mathcal I}^3 \oplus {\mathcal I}W/{\mathcal I}^2 W \otimes_X  {\mathcal I}/{\mathcal I}^2 \rar \frac{ W/{\mathcal I}^2 W \otimes_X  {\mathcal I}/{\mathcal I}^3}{{\mathcal I} W/{\mathcal I}^2 W \otimes_X  {\mathcal I}^2/{\mathcal I}^3} \rar W/{\mathcal I} W \otimes_X  {\mathcal I}/{\mathcal I}^2 \rar 0$. Quotienting out by the kernel of the map $\id +  p_2$ to $W/{\mathcal I} W \otimes_X  {\mathcal I}^2/{\mathcal I}^3$ which we demote by $ker$ we see that $\id_V \otimes_X \psi_{X/S} + p_2 \circ (\theta_V \otimes \id_{\Omega_{X/S}})$ is given by the exact sequence $0 \rar W/{\mathcal I} W \otimes_X  {\mathcal I}^2/{\mathcal I}^3 \rar \frac{ W/{\mathcal I}^2 W \otimes_X  {\mathcal I}/{\mathcal I}^3}{{\mathcal I} W/{\mathcal I}^2 W \otimes_X  {\mathcal I}^2/{\mathcal I}^3}/ker \rar  W/{\mathcal I} W \otimes_X  {\mathcal I}/{\mathcal I}^2 \rar 0$. The desired result follows from the fact that we have a surjection $\varphi:F:= \frac{ W/{\mathcal I}^2 W \otimes_X  {\mathcal I}/{\mathcal I}^3}{{\mathcal I} W/{\mathcal I}^2 W \otimes_X  {\mathcal I}^2/{\mathcal I}^3}/ker \rar  W/{\mathcal I}^2 W \otimes {\mathcal I}/{\mathcal I}^3 =  W/{\mathcal I}^3 W$. This gives rise to the following commutative diagram:\\

\[
 \begin{CD}
   0 @> >>  W/{\mathcal I} W \otimes_X  {\mathcal I}^2/{\mathcal I}^3 @> >> F @> >> W/{\mathcal I} W \otimes_X  {\mathcal I}/{\mathcal I}^2 @> >> 0\\
@. @VV{\id}V @VV{\varphi}V @V{\id}VV @.\\
 0 @> >>  W/{\mathcal I} W \otimes_X  {\mathcal I}^2/{\mathcal I}^3 @> >>  W/{\mathcal I}^3 W @> >> W/{\mathcal I} W \otimes_X  {\mathcal I}/{\mathcal I}^2 @> >> 0\\
\end{CD}
 \]\\
\end{proof}

\begin{claim}
$\tilde{\ttt_2}(V) = \tilde{\ch_2}(V) - (\id_V \otimes \psi_{X/S}) \circ \theta_V$
\end{claim}
\begin{proof}
 $\tilde{\ttt_2}(V) = \tilde{\ch_2}(V) + p_2 \circ \tilde{\ttt_2}(V)$.But $\tilde{\ttt_2}(V) = (\theta_V \otimes \id_{\Omega}) \circ theta_V$ and $\psi_V \circ \theta_V =0$. Thus, $\psi_V \circ \theta_V = (\id_V \otimes \psi_{X/S}) \circ \theta_V + p_2 \circ (\theta_V \otimes \id_{\Omega}) \circ \theta_V =0 \implies  p_2 \circ \tilde{\ttt_2}(V)= - (\id_V \otimes \psi_{X/S}) \circ \theta_V$.\\
\end{proof}

We now note that $\tilde{\ttt_3}(V) = (\tilde{\ttt_2}(V) \otimes \id_{\Omega}) \circ \theta_V = (\tilde{\ch_2}(V) \otimes \id_{\Omega}) \circ \theta_V -((\id_V \otimes \psi_{X/S}) \otimes \id_{\Omega}) \circ (\theta_V \otimes \id_{\Omega}) \circ theta_V = (\tilde{\ch_2}(V) \otimes \id_{\Omega}) \circ \theta_V -((\id_V \otimes \psi_{X/S}) \otimes \id_{\Omega}) \circ (\tilde{\ch_2}(V) - (\id_V \otimes \psi_{X/S}) \circ \theta_V)$. Also, as $\ttt_3(V) \in \Hm^3(X , {\Omega}^{\otimes 3})$ actually lies in $\Hm^3(X , {{\Omega}^{\otimes 3}}_{\tau})$ where ${{\Omega}^{\otimes 3}}_{\tau}$ denotes the part of ${\Omega}^{\otimes 3}$ that is invariant under the action of the $3$-cycle, i.e, $\sym^3 \Omega \oplus {\wedge}^3 \Omega$, we only need to see what part of the right hand side of the above equation for $\tilde{\ttt_3}(V)$ is left when we project to  $\sym^3 \Omega \oplus {\wedge}^3 \Omega$. $((\id_V \otimes \psi_{X/S}) \otimes \id_{\Omega}) \circ (\tilde{\ch_2}(V) - (\id_V \otimes \psi_{X/S}) \circ \theta_V) \in \eee^3(V,V \otimes \sym^2 \Omega \otimes \Omega)$ , which means that its projection to $\eee^3(V, V \otimes {\wedge}^3 \Omega)$ vanishes. Also, $(\tilde{\ch_2}(V) \otimes \id_{\Omega}) \circ \theta_V \in \eee^3(V,V \otimes {\wedge}^2 \Omega \otimes \Omega)$ which means that its projection to $\eee^3(V, V \otimes \sym^3 \Omega)$ vanishes. From this, we see that upto expressions that vanish on taking the trace, $\tilde{\ttt_3}(V) = \tilde{\ch_3}(V) - p_3 \circ ((\id_V \otimes \psi_{X/S}) \otimes \id_{\Omega}) \circ (\tilde{\ch_2}(V) - (\id_V \otimes \psi_{X/S}) \circ \theta_V)$.It remains to be shown that  $ p_3 \circ ((\id_V \otimes \psi_{X/S}) \otimes \id_{\Omega}) \circ  (\id_V \otimes \psi_{X/S}) \circ \theta_V =0$. This follows from the following claim: \\

\begin{claim}
$ p_3 \circ (( \psi_{X/S} \otimes \id_{\Omega}) \circ   \psi_{X/S})$
\end{claim}
\begin{proof}
Note that $\psi_{X/S}$ is given by the exact sequence $0 \rar {\mathcal I}^2/{\mathcal I}^3 \rar {\mathcal I}/{\mathcal I}^3 \rar {\mathcal I}/{\mathcal I}^2 \rar 0$.We therefore get the following commutative diagram whose top row is the exact sequence giving $\psi_{X/S} \otimes \id_{\Omega}$ and whose bottom row is the exact sequence giving  $p_3 \circ  (\psi_{X/S} \otimes \id_{\Omega})$: \\

\[
\begin{CD}
0 @> >> {\mathcal I}/{\mathcal I}^2 \otimes {\mathcal I}^2/{\mathcal I}^3  @> i >> {\mathcal I}/{\mathcal I}^2 \otimes {\mathcal I}/{\mathcal I}^3 @> g >> {\mathcal I}/{\mathcal I}^2 \otimes {\mathcal I}/{\mathcal I}^2 @> >>0 \\
@. @AA{\id}A @AAA @AAA @. \\
0 @> >> {\mathcal I}/{\mathcal I}^2 \otimes {\mathcal I}^2/{\mathcal I}^3 @> >> W_1 @> k >> {\mathcal I}^2/{\mathcal I}^3 @> >> 0\\
@. @VVhV @VVV @V{\id}VV @.\\
0 @> >> {\mathcal I}^3/{\mathcal I}^4 @> >> W_2 @> >> {\mathcal I}^2/{\mathcal I}^3 @> >> 0\\
\end{CD}
\]\\

Therefore, $p_3 \circ  (\psi_{X/S} \otimes \id_{\Omega}) \circ \psi_{X/S}$ is given by the exact sequence $0 \rar {\mathcal I}^3/{\mathcal I}^4 \rar W_2  \rar {\mathcal I}/{\mathcal I}^3 \rar {\mathcal I}/{\mathcal I}^2 \rar 0$. Suppose we show that it is also given by the exact sequence $0 \rar {\mathcal I}^3/{\mathcal I}^4 \rar {\mathcal I}^2/{\mathcal I}^4  \rar {\mathcal I}/{\mathcal I}^3 \rar {\mathcal I}/{\mathcal I}^2 \rar 0$ , then it is $0$ since the following diagram commutes:\\

\[
 \begin{CD}
  0 @> >> {\mathcal I}^3/{\mathcal I}^4  @> >> {\mathcal I}^2/{\mathcal I}^4 @> >> {\mathcal I}/{\mathcal I}^3 @> >> {\mathcal I}/{\mathcal I}^2 @> >>0\\
@. @AAA @AAA @AAA @AAA @.\\
0 @> >> 0 @> >> {\mathcal I}^2/{\mathcal I}^4  @> >> {\mathcal I}/{\mathcal I}^4 @ > >> {\mathcal I}/{\mathcal I}^2 @> >>0\\
\end{CD}
\]\\

It is therefore, enough to find an arrow $f$ so that the following diagram commutes:\\

\[
 \begin{CD}
0 @> >>  {\mathcal I}^3/{\mathcal I}^4  @> >> W_2 @> >> {\mathcal I}^2/{\mathcal I}^3 @> >> 0\\
@. @VV{\id}V @VVfV @V{\id}VV @.\\
  0 @> >>  {\mathcal I}^3/{\mathcal I}^4  @> >>  {\mathcal I}^2/{\mathcal I}^4  @> >> {\mathcal I}^2/{\mathcal I}^3 @> >> 0\\
  \end{CD}
\]\\

Here, $W_2 = \frac{ g^{-1}(\sym^2  {\mathcal I}/{\mathcal I}^2) \oplus  {\mathcal I}^3/{\mathcal I}^4}{(i,-h)( {\mathcal I}^2/{\mathcal I}^3 \otimes  {\mathcal I}/{\mathcal I}^2)}$. We therefore, need to be able to find a map from $g^{-1}(\sym^2  {\mathcal I}/{\mathcal I}^2)$ to ${\mathcal I}^2/{\mathcal I}^4$. This we will do locally. Observe that the morphism $\sym^2 {\mathcal I}/{\mathcal I}^2 \rar {\mathcal I}^2/{\mathcal I}^3$ factors as $\sym^2 {\mathcal I}/{\mathcal I}^2 \rar {\mathcal I}/{\mathcal I}^2 \otimes {\mathcal I}/{\mathcal I}^2 \rar {\mathcal I}^2/{\mathcal I}^3$. So we have a map $g^{-1}(\sym^2 {\mathcal I}/{\mathcal I}^2) \rar  {\mathcal I}/{\mathcal I}^3 \otimes {\mathcal I}/{\mathcal I}^2 \rar {\mathcal I}^2/{\mathcal I}^3$. The kernel of the composite is precisely ${\mathcal I}^2/{\mathcal I}^3 \otimes {\mathcal I}/{\mathcal I}^2$ and it factors through ${\mathcal I}^2/{\mathcal I}^4$. The rest is easy.

\end{proof}
\end{proof}

\subsection{Computation of $\psi_{G/S}$ where $G=G(r,n)$}

We revert now to the notation of Section 4. By Claim 3 of the previous subsection , $ p_2 \circ \tilde{\ttt_2}(Q) = -(\id_Q \otimes \psi_{G/S}) \circ \theta_Q$. As an element of $\hhh_K( {\wedge}^2 \Omega \otimes Q, {\Omega}^2 \otimes Q)$, we have $$-(\id_Q \otimes \psi_{G/S}) \circ \theta_Q = p_2 \circ \sum_{l_1,l_2;m_1,m_2;r_1,r_2} (f_{m_1} \otimes e_{l_1}) \circ (f_{m_2} \otimes e_{l_2}) \bigotimes (e_{m_1} \otimes v_{r_1}) \otimes (e_{m_2} \otimes v_{r_2}) \bigotimes (f_{l_1} \otimes u_{r_1}) \otimes (f_{l_2} \otimes u_{r_2})$$ $$ = \frac{1}{2} \sum_{l_1,l_2;m_1,m_2;r_1,r_2} (f_{m_1} \otimes e_{l_1}) \circ (f_{m_2} \otimes e_{l_2}) \bigotimes (e_{m_1} \otimes v_{r_1}) \otimes (e_{m_2} \otimes v_{r_2})$$ $$ \bigotimes [(f_{l_1} \otimes u_{r_1}) \otimes (f_{l_2} \otimes u_{r_2}) + (f_{l_2} \otimes u_{r_2}) \otimes (f_{l_1} \otimes u_{r_1})]$$ \\
 Also, $\theta_Q = \sum_{l_1,m_1,r_1} (f_{m_1} \otimes e_{l_1}) \bigotimes (e_{m_1} \otimes v_{r_1}) \bigotimes (f_{l_1} \otimes u_{r_1})$ in $\hhh_K( Q \otimes \Omega, Q \otimes \Omega)$. Therefore,if $\psi_{G/S}$ is given by an element $\varphi \in \hhh_K( \Omega \otimes \Omega, \Omega \otimes \Omega)$, we observe that the following diagrams commute:\\
\[
 \begin{CD}
0 @> >>Q \otimes {\Omega}^{\otimes 2} @> >> Q \otimes E @> >> Q \otimes {\Omega} @> >> 0 \\
@. @AA(\bar{\id_Q \otimes \varphi})A @AAA @A{\id}AA \\
...@> >> Q \otimes {\Omega}^{\otimes 2} \otimes \sym^* \Omega @> \alpha_1 >> Q \otimes {\Omega} \otimes \sym^* \Omega @> >> Q \otimes {\Omega} @> >> 0 \\
\end{CD}
\]\\

\[
\begin{CD}
  0 @> >>Q \otimes {\Omega} @> >> J_1(Q) @> >> Q  @> >> 0 \\
 @. @AA{\theta_Q}A @AAA @A{\id}AA @.\\
...@> >> Q \otimes {\Omega} \otimes \sym^* \Omega @> \alpha_1 >> Q \otimes  \sym^* \Omega @> >> Q  @> >> 0 \\
\end{CD}
\]\\

The following diagran, therefore, commutes:\\

\[
 \begin{CD}
0 @> >>Q \otimes {\Omega}^{\otimes 2} @> >> Q \otimes E @> >> Q \otimes {\Omega} @> >> 0 \\
@. @AA(\id_Q \otimes \varphi)A @AA{\gamma}A @A{\id}AA @.\\
...@> >> Q \otimes {\Omega}^{\otimes 2} \otimes \sym^* \Omega @> \alpha_1 >> Q \otimes {\Omega} \otimes \sym^* \Omega @> >> Q \otimes {\Omega} @> >> 0 \\
@. @AA{\id_{\Omega} \otimes \tilde{\theta}_Q}A @AA{\tilde{\theta}_Q}A @A{\theta_Q}AA @.\\
...@> >> Q \otimes {\Omega}^{\otimes 2} \otimes \sym^* \Omega @> \alpha_1 >> Q \otimes {\Omega} \otimes \sym^* \Omega @> >> Q \otimes {\Omega} @> >> 0 \\
\end{CD}
\]\\

This gives rise to the following commuting diagram:\\

\[
\begin{CD}
 0 @> >>Q \otimes {\Omega}^{\otimes 2} @> >> Q \otimes E @> >>J_1(Q) @> >> 0\\
@. @AA{ (\id_Q \otimes \varphi) \circ (\id_{\Omega} \otimes \theta_Q)}A @AA{ \gamma \circ \theta_Q}A @A{\beta}AA @.\\
...@> >> Q \otimes {\Omega}^{\otimes 2} \otimes \sym^* \Omega @> \alpha_1 >> Q \otimes {\Omega} \otimes \sym^* \Omega @> >> Q \otimes \sym^* {\Omega} @> >> 0 \\
\end{CD}
\]\\

Therefore, $ (\id_Q \otimes \varphi) \circ (\id_{\Omega} \otimes \theta_Q) \circ i = p_2 \circ \tilde{\ttt_2}(Q)$ as an element of $\hhh_K( {\wedge}^2 \Omega, {\Omega}^{\otimes 2})$. We can use this to determine $\varphi$. Note what $ p_2 \circ \tilde{\ttt_2}(Q)$ does to a basis element $e_a \bigotimes (f_b \otimes u_c) \bigotimes (f_{b'} \otimes u_{c'})$ : $$p_2 \circ \tilde{\ttt_2}(Q) (e_a \bigotimes (f_b \otimes u_c) \bigotimes (f_{b'} \otimes u_{c'}))$$ $$ = \sum_{l_1,l_2;m_1,m_2;r_1,r_2} \delta_{am_2}\delta_{l_2m_1}\delta_{bm_1}\delta_{m_2b'}\delta_{cr_1}\delta_{r_2c'} e_{l_1} \bigotimes \frac{f_{l_1} \otimes u_{r_1} \bigotimes f_{l_2} \otimes u_{r_2} + f_{l_2} \otimes u_{r_2} \bigotimes f_{l_1} \otimes u_{r_1}}{2}$$ A term here is nonzero if $m_2 = a = b'$, $m_1 = b= l_2$,$r_1 = c$ and $r_2 = c'$. Therefore, $$p_2 \circ \tilde{\ttt_2}(Q) (e_a \bigotimes (f_b \otimes u_c) \bigotimes (f_{b'} \otimes u_{c'}))$$ $$ = \frac{\delta_{ab'}}{2} \sum_{l_1} e_{l_1} \bigotimes  (f_{l_1} \otimes u_{c} \bigotimes f_{b} \otimes u_{c'} + f_{b} \otimes u_{c'} \bigotimes f_{l_1} \otimes u_{c})$$  On the other hand, we know that $$\id_{\Omega} \otimes \theta_Q ( e_a \bigotimes (f_b \otimes u_c) \bigotimes (f_{b'} \otimes u_{c'}))$$ $$ = \sum_{l_1,m_1,r_1} \delta_{am_1}\delta_{b'm_1}\delta_{c'r_1} e_{l_1} \bigotimes (f_b \otimes u_c) \bigotimes (f_{m_1} \otimes u_{r_1})$$ $$ = \delta_{ab'} \sum_{l_1} e_{l_1} \bigotimes (f_b \otimes u_c) \bigotimes (f_{l_1} \otimes u_{c'}) $$  Therefore, if $b \neq a$, $$p_2 \circ \tilde{\ttt_2}(Q) (e_a \bigotimes (f_b \otimes u_c) \bigotimes (f_a \otimes u_{c'}) - e_a \bigotimes (f_a \otimes u_{c'}) \bigotimes (f_b \otimes u_{c}))$$ $$ = \frac{1}{2} \sum_{l_1}  e_{l_1} \bigotimes [(f_{l_1} \otimes u_c) \bigotimes (f_b \otimes u_{c'} + (f_b \otimes u_{c'}) \bigotimes (f_{l_1} \otimes u_{c}] $$ $$ = (\id_Q \otimes \varphi) \circ (\id_{\Omega} \otimes \theta_Q) (e_a \bigotimes (f_b \otimes u_c) \bigotimes (f_a \otimes u_{c'}) - e_a \bigotimes (f_a \otimes u_{c'}) \bigotimes (f_b \otimes u_{c}))$$ $$   = (\id_Q \otimes \varphi) \sum_{l_1} e_{l_1} \bigotimes (f_b \otimes u_c) \bigotimes (f_{l_1} \otimes u_{c'})$$ $$ \implies \varphi ( (f_b \otimes u_c) \bigotimes (f_{l_1} \otimes u_{c'}))$$ $$ = \frac{1}{2} [(f_{l_1} \otimes u_c) \bigotimes (f_b \otimes u_{c'} + (f_b \otimes u_{c'}) \bigotimes (f_{l_1} \otimes u_{c})]$$  Therefore, as an element in $\enn_K({Q^*}^{\otimes 2} \otimes S^{\otimes 2})$, we see that $\psi_{G/S}$ is given by $ {(1 \text{ } 2)}_{Q^*} \otimes \id_S + \id_{Q^*} \otimes {(1 \text{ } 2)}_{S}$. \\

\subsection{Formula for $\ttt_k(Q)$ in terms of $\ch_l(Q)$ , $l \leq k$}

First, we compute $({\id_{\Omega}}^{\otimes r} \otimes \psi_{G/S} \otimes {\id_{\Omega}}^{\otimes k-2-r}) \circ \ttt_{k-1}(Q)$. Note that if $\psi_{G/S}$ is given by the exact sequence $ 0 \rar {\Omega}^{\otimes 2} \rar E \rar \Omega \rar 0$, then ${\id_{\Omega}}^{\otimes r} \otimes \psi_{G/S} \otimes {\id_{\Omega}}^{\otimes k-2-r}$ is given by $  0 \rar {\Omega}^{\otimes k} \rar {\Omega}^{\otimes r} \otimes E \otimes {\Omega}^{\otimes k-2-r} \rar {\Omega}^{\otimes k-1} \rar 0$. By the computation of the previous subsection we obtain the following commutative diagram: \\

\[
 \begin{CD}
    0 @> >>  {\Omega}^{\otimes k} @> >> {\Omega}^{\otimes r} \otimes E \otimes {\Omega}^{\otimes k-2-r} @> >> {\Omega}^{\otimes k-1} @> >> 0\\
@. @AA{\varphi_r}A @AAA @A{\id}AA @.\\
...@> >> \Omega^{\otimes k} \otimes \sym^* \Omega @> \alpha_{r+1} >> \Omega^{\otimes k-1} \otimes \sym^* \Omega @> >>  \Omega^{\otimes k-1} @> >> 0\\
\end{CD}
\]\\

By Observation 1 of Section 4.3 , ${(-1)}^r \alpha_{r+1}$ commutes with the Koszul differential. We can therefore use the following diagram for our computation:\\

\[
 \begin{CD}
    0 @> >>  {\Omega}^{\otimes k} @> >> {\Omega}^{\otimes r} \otimes E \otimes {\Omega}^{\otimes k-2-r} @> >> {\Omega}^{\otimes k-1} @> >> 0\\
@. @AA{{(-1)}^r \varphi_r}A @AAA @A{\id}AA @.\\
...@> >> \Omega^{\otimes k} \otimes \sym^* \Omega @> { {(-1)}^r \alpha_{r+1}} >> \Omega^{\otimes k-1} \otimes \sym^* \Omega @> >>  \Omega^{\otimes k-1} @> >> 0\\
\end{CD}
\]\\

By the previous subsection, $\varphi_r = {(r+1 \text{ } r+2)}_{Q^*} \otimes \id_S + \id_{Q^*} \otimes {(r+1 \text{ } r+2)}_S$ as an element of $\enn_K({Q^*}^{\otimes k} \otimes S^{\otimes k})$. Following the method of computation of Section 4, it is clear that  $({\id_{\Omega}}^{\otimes r} \otimes \psi_{G/S} \otimes {\id_{\Omega}}^{\otimes k-2-r}) \circ \ttt_{k-1}(Q) = \sum_{\sigma \in S_k} \sn(\sigma) (\sigma \otimes \sigma) (\tau_{k-1} \otimes \id) ((r+1 \text{ } r+1) \otimes \id + \id \otimes (r+1 \text{ } r+2)) {(-1)}^r$ as an element of $\enn_K({Q^*}^{\otimes k} \otimes S^{\otimes k})$. Again, it can be seen easily that this also equals ${(-1)}^r \sum_{\sigma \in S_k} \sn(\sigma) (\sigma \otimes \sigma) ([\tau_{k-1} (r+1 \text{ } r+2) - (r+1 \text{ } r+2) \tau_{k-1}] \otimes \id)$. But $\tau_{k-1} (r+1 \text{ } r+2) = \tau_k$. Therefore,  $({\id_{\Omega}}^{\otimes r} \otimes \psi_{G/S} \otimes {\id_{\Omega}}^{\otimes k-2-r}) \circ \ttt_{k-1}(Q) = {(-1)}^r \sum_{\sigma \in S_k} \sn(\sigma) (\sigma \otimes \sigma) [(\tau_k - (r+1 \text{ } r+2) \tau_k (r+1 \text{ } r+2) ) \otimes \id]$. This is true for all $r \leq k-2$. Observe that $\ttt_k(Q) - \ch_k(Q) = C \sum_{\sigma \in S_k} \sn(\sigma) (\sigma \otimes \sigma) [\sum_{\omega \in S_k} (\tau_k - \omega \tau_k {\omega}^{-1}) \otimes \id]$, where $C$ is a constant. The transpositions $\mu_i = (i \text{ } i+1)$ , $1 \leq i \leq k-1$ generate $S_k$. Also, following the proof of Lemma 14 we see that $\omega_* \sum_{\sigma \in S_k} \sn(\sigma) (\sigma \otimes \sigma) ( \gamma \otimes \id) = \sn(\omega) \sum_{\sigma \in S_k} \sn(\sigma) (\sigma \otimes \sigma) ( \omega^{-1} \gamma \omega \otimes \id)$. It follows that if $\omega = \mu_{i_1}.... \mu_{i_s}$, since $\tau_k - \omega \tau_k \omega^{-1} = \tau_k - \mu_{i_1} \tau_k \mu_{i_1} + mu_{i_1} ( \tau_k - \lambda \tau_k \lambda^{-1}) \mu_{i_1}$, where $\lambda = \mu_{i_2} ... \mu_{i_s}$, we can write $\tau_k - \omega \tau_k \omega^{-1} = \tau_k - \mu_{i_1} \tau_k \mu_{i_1} + \beta_{2}*(\tau_k - \mu_{i_2}*\tau_k) + ... + \beta_{s}*(\tau_k - \mu_{i_s}*\tau_k)$, where $\beta_i \in S_k$ for all $i$and $*$ denotes conjugation. Therefore, $\ttt_k(Q) -\ch_k(Q) = \sum \gamma_{j*} A_j \circ \ttt_{k-1}(Q)$ where $\gamma_j \in KS_k$ and $A_j =  ({\id_{\Omega}}^{\otimes j-1} \otimes \psi_{G/S} \otimes {\id_{\Omega}}^{\otimes k-1-j}), \text{ } 1 \leq j \leq k-1$. By induction, there exists a formula for $\ttt_k(Q)$ in terms of $\ch_l(Q), \text{ } l \leq k$ and $\psi_{G/S}$. \\

\subsection{Proving the formula for $\ttt_k(V)$}

Let $X$ be a projective variety. Let $A_r$ denote the functor $({\id_{\Omega}}^{\otimes r-1} \otimes \psi \otimes {\id_{\Omega}}^{\otimes k-1-r})$. We need to check that $A_r \circ \ttt_{k-1}$ is functorial with respect to pull-backs for $1 \leq r \leq k-1$. For this, we note that the diagram below commutes given a morphism $f:X \rar Y$ where $X$ and $Y$ are projective varieties over $S$:\\

\[
\begin{CD}
 0 @> >> {f^* \Omega_Y}^{\otimes k} @> >> {f^* \Omega_Y}^{\otimes r-1 } \otimes f^* E_Y \otimes {f^* \Omega_Y}^{\otimes k-1-r} @> >> ...@> >> \calo_X @> >>0\\
@. @AAA @AAA @AAA @A{\id}AA @. \\
0 @> >> {\Omega_X}^{\otimes k} @> >> { \Omega_X}^{\otimes r-1 } \otimes  E_X \otimes { \Omega_X}^{\otimes k-1-r} @> >> ...@> >> \calo_X @> >>0\\
\end{CD}
\]\\

It follows that the exact sequence on the bottom row is the element of $\eee^k(\calo_X , {\Omega_X}^{\otimes k})$ induced from the element of $\eee^k(\calo_X,{f^* \Omega_Y}^{\otimes k})$ given by the top row under the map ${f^* \Omega_Y}^{\otimes k} \rar {\Omega_X}^{\otimes k}$. By the definition of $f^*:\eee^k(\calo_Y , {\Omega_Y}^{\otimes k}) \rar \eee^k(\calo_X , {\Omega_X}^{\otimes k})$, it follows that $A_r \circ \ttt_{k-1}$ is functorial under pullbacks. \\
In the previous subsection, we showed that $\ttt_k(Q) = \ch_k(Q) + \sum_{1 \leq j \leq k-1} \gamma_{j*} A_j \circ \ttt_{k-1}(Q)$ where $\gamma_j \in KS_k , \text{ } 1 \leq j \leq k-1$. As $\gamma_{j*} A_j \circ \ttt_{k-1}$ is functorial under pullbacks, if $V$ is a vector bundle on $X$ a projective variety, we see that for all sufficiently large $n$, $\ttt_k(V \otimes \calo(n)) = \ch_k(V \otimes \calo(n)) + \sum_{1 \leq j \leq k-1} \gamma_{j*} A_j \circ \ttt_{k-1}(V \otimes \calo(n))$. Both sides of this equation are polynomials in $n$ whose constant terms are $\ttt_k(V)$ and $ \ch_k(V) + \sum_{1 \leq j \leq k-1} \gamma_{j*} A_j \circ \ttt_{k-1}(V)$ respectively. It follows that $\ttt_k(V) = \ch_k(V) + \sum_{1 \leq j \leq k-1} \gamma_{j*} A_j \circ \ttt_{k-1}(V)$. Proceeding inductively as in the previous section, we obtain a formula for $\ttt_k(V)$ in terms of $\ch_l(V) \text{ } l \leq k$ and $\psi_{X/S}$. We have thus, proven the following theorem:

\begin{theorem}
If $X$ is a projective variety and $V$ is a vector bundle on $X$, then $\ttt_k(V) = \ch_K(V) + \sum_{l \leq k} D_{kl} \circ \ch_l(V)$ where $D_{kl}$ are elements of $\eee^{k-l}(\Omega^{\otimes l} ,\Omega^{\otimes k})$ which are functorial under pullbacks. Moreover $D_{kl}$ is obtained from $\psi_{X/S}$ by applying finitely many of the following operations:\\
1: Tensoring with $\id_{\Omega}$\\
2: Yoneda multiplication\\
3: Action of permutation group elements on $\eee^j(-, \Omega^{\otimes j})$.
\end{theorem}

\subsection{Proper subfunctors of the Hodge functors $\Hm^q(X, \Omega^p)$, $p,q \geq 2$}

The formula for $\ttt_k$ in terms of $\ch_l, \text{ } l \leq k$ also easily gives us a method for finding an increasing chain of proper contrvariant subfunctors of the Hodge functors $\Hm^q(X, \Omega^p)$ for projective varieties over field of characteristic $0$.  Note that $D_{kl}$ maps $\Hm^l(X, \Omega^l)$ to $\Hm^k(X, \Omega^{\otimes k})$, and, by the formula for $\ttt_k$ given in the previous subsection, $D_{kl} \circ \alpha_l(V) = \ttt_k(\alpha_l(V))$. If $X = G(r,n)$ a Grassmannian, and $V = Q$ , the universal quotient bundle of $G(r,n)$ and $n$ is large enough, we see that $D_{kl} \circ \alpha_l(Q) = \ttt_k( \alpha_l(Q)) \neq 0$ if $l \geq 2$. Therefore, Yoneda composition with $D_{kl}$ does not kill $\Hm^{k,k}$ in general. On the other hand, it was proven in Nori[1] that $\psi_{X/S} = 0$ if $X$ is an Abelian variety, (a torus for example). Therefore, if $X= G(r,n) \times T$ where $T$ is a torus, then $\alpha_l(p_1^*Q)$ is not in the kernel of $D_{kl} \circ $, but $D_{kl} \circ p_2^* Y = 0 \text{ } \forall Y \in \Hm^{l,l}(T) $. Therefore, we see that $\Hm^{l,l}_k := \ker D_{kl} \circ : \Hm^{l,l} \rar \Hm^k(X, \Omega^{\otimes k})$ is a proper subfunctor of $\Hm^{l,l}$ for all $l \geq 2$. \\
Similarly, we see that $\Hm^{p,q}_k$ given by $\ker (D_{kq} \otimes \id_{\Omega}^{\otimes p-q}) \circ : \Hm^{p,q} \rar  \Hm^k(X, \Omega^{\otimes k+p-q})$ if $p > q$ and  $\ker D_{kp} \circ : \Hm^{p,q} \rar  \Hm^{k+q-p}(X, \Omega^{\otimes k})$ otherwise is a proper subfunctor of $\Hm^{p,q}$. To see this, again consider the case when $X = G(r,n) \times T$ as before, $T$ a suitable torus. If $p > q$ consder the element $ \alpha_q(p_1^*Q) \cup p_2^*Y \in \Hm^{p,q}(X)$ where $Y$ is a nonzero element of $\Hm^{p-q,0}(T)$. One checks that $(D_{kq} \otimes \id_{\Omega}^{\otimes p-q}) \circ (\alpha_q(p_1^*Q) \cup p_2^*Y) = (D_{kq} \circ \alpha_q(p_1^*Q)) \cup  p_2^*Y = \ttt_k(\alpha_q(Q)) \cup p_2^*Y \neq 0$. On the other hand, if $Z \in \Hm^{p,q}(T)$, then  $(D_{kq} \otimes \id_{\Omega}^{\otimes p-q}) \circ p_2^*Z =0$. This proves that $\Hm^{p,q}_k$ is a proper subfunctor of $\Hm^{p.q}$ if $p>q$. If $p<q$ we note that if $Y \in \Hm^{0,q-p}(T)$ is nonzero, then $D_{kp} \circ ( \alpha_p(p_1^*Q) \cup p_2^*Y) = (D_{kp} \circ  \alpha_p(p_1^*Q)) \cup p_2^*Y = \ttt_k(\alpha_p(p_1^*Q)) \cup p_2^*Y \neq 0$ and that if $Z \in \Hm^{p,q}(T)$ then  $D_{kp} \circ p_2^*Z =0$. This proves that $\Hm^{p,q}_k$ is a proper subfunctor of $\Hm^{p,q}$ for all $k > q$ where $p,q \geq 2$. \\
In fact, a little more is true. We mentioned in the previous section that $\ttt_k(V) = \ch_k(V) + \sum_{1 \leq j \leq k-1} \gamma_{j*} A_j \circ \ttt_{k-1}(V)$. It follows that $D_{kl} = \sum_{1 \leq j \leq k-1} \gamma_{j*} A_j \circ D_{k-1 l} $ if $l < k-1$. This in fact tells us that the subfunctors $\Hm^{p,q}_k$ in fact form an increasing chain of proper subfunctors of $\Hm^{p,q}$ for $k > \min{\{p,q\}}$ (as a theory).

\subsection{future problems}

Some questions arising out of this work that need to be addressed in the future are described below.\\

1. The filtration $F_r \CH^l$ described in Theorem 4 was shown to be nontrivial
as a theory. However, on going through the proof, one gets a feeling that more can be shown. Intuitively,I feel that it is possible to show that this filtration is strictly increasing as a theory and that given any $l \geq 2$ fixed, and $r \geq 2$ , there exists some Grassmannian $G=G(r,n)$ so that $\alpha_l(Q) \in F_r \CH^l(G) \setminus F_{r-1} \CH^l(G)$. One approach to this question is entirely combinatorial and boils down to showing that for some $k$ and a particular $\beta \in KS_k$ depending on $l$ and $k$ only, the subspace spanned by the conjugates of $\beta_{r-1}$ is of strictly smaller dimension than that spanned by conjugates of $\beta_r$. Here, $\beta_i$ is the image of $\beta$ under the projection $KS_k \rar \oplus_{|\lambda| \leq i} End(V_{\lambda})$. Approaching this question along these lines would indeed involve algebraic combinatorics extensively.    \\

2. Is the map $\D_{kl} \circ - : \Hm^0(X, \Omega^l) \rightarrow \Hm^{k-l}(X, \Omega^{\otimes k})$ nonzero ? If it is zero, then are there other elements $\D_m$ in $\eee^{m-l}(\Omega^{\otimes l}, \Omega^{\otimes m})$ which are obtained from $\psi_{X/S}$ by tensoring with $\id_{\Omega}$, Yoneda multiplication and the action of permutation group elements on $\eee^j(-, \Omega^{\otimes j})$ in finitely many steps, so that $\D_m \circ -: \Hm^0(X, \Omega^l) \rightarrow \Hm^m(X, \Omega^{\otimes k})$ nonzero ? \\
 When $l=1$, the answer is in the negative. This follows from the fact that the map $X \rightarrow J_X$ induces an isomorphism on $\Hm^{1,0}$ where $J_X$ is the Albanese variety of $X$ and the fact that $\psi_{A/S} = 0$ if $A$ is an Abelian variety. However, the above question is open for $l \geq 2$. One can also see that the map $D_{kl} \circ -$ is not injective in general. For example, it is zero on any Abelian variety, and hence, on elements of $\Hm^{l,0}$ that are pullbacks of such elements living on $\Hm^{1,0}$ of an Abelian variety. \\
Proving that these maps are nonzero, will therefore give us subfunctors of the $\Hm^{l,0}$. At least for the case when $k=3, l=2$, we could obtain elements of $\Hm^{2,0}$ of some surfaces explicitly, Yoneda compose them with $D_{32}$ and see if the resulting short exact sequence splits.\\

3. The map $\D_{k1} \circ -: \Hm^{1,1} \rightarrow  \Hm^{k-1}(X, \Omega^{\otimes k})$ is known to be zero on the span of the image of $\ch_1$. What happens to this map on the rest of $\Hm^{1,1}$ ?. \\

4. The functoriality of the $\Hm^{p,q}_k$ introduce new restrictions the pullback maps on cohomology associated with morphisms between projective varieties must satisfy. Identifying these restrictions concretely might enable us to obtain more results regarding the nonexistence of morphisms. The problem of computing the invariants $\hm^{l,l}_k := \dim \Hm^{l,l}_k$ for Grassmannians looks like a combinatorial problem at present, though other approaches may be needed in the future. If it is not possible to find these invariants exactly, estimating them may also be useful in obtaining more nonexistence of morphisms results. In addition, if question 1. is answered in the affirmative, the same problem arises regarding the new invariants that that answer will lead us to. \\


\begin{thebibliography}{GGG}

\bibitem{1} Nori, M. V. {\em Personal Communication}

\bibitem{2} Loday, J- L. {\em Cyclic Homology}, Springer ,1997.

\bibitem{3} Fulton, W. and Harris, J. {\em Representation Theory: a first course} , Springer Verlag, GTM 129, 1991.

\bibitem{4} Bott, Raoul. {\em Homogenous Vector Bundles}, Annals of Mathematics, 2nd Ser., Vol 66,No.2 (Sep., 1957), pp 203-248.

\bibitem{5} Reutenauer, Christophe. {\em Free Lie Algebras}, London MathematicalSociety Monographs, new series, no. 7, Oxford University Press, 1993.

\bibitem{6} Kapranov, M. Rozansky-Witten invariants via Atiyah classes. {\em Compositio Math.} 115(1999), 71-113.

\bibitem{7} Paranjape, K. H.; Srinivas, V.  Self maps of homogeneous spaces . {\em Invent. Math.} { 98}(1989), no.2,425-444

\bibitem{8}  Paranjape, K. H.; Srinivas, V. Continuous self-maps of quadric hypersurfaces. {\em Proceedings of the Indo-French conference on Geometry} ({\em Bombay}, 1989), 135-148, {\em Hindustan Book Agency, Delhi},1993.

\end{thebibliography}
\end{document}